%% file: header.tex
\documentclass[12pt]{article}

\usepackage{latexsym,a4,AM,epsfig}

\def\Bbb{\bf}

\newcommand{\E}{e}

\newtheorem{aslemma}[subsection]{Lemma}{\bf}{\it}
\newtheorem{astheorem}[subsection]{Theorem}{\bf}{\it}
\newtheorem{asprop}[subsection]{Proposition}{\bf}{\it}
\newtheorem{ascor}[subsection]{Corollary}{\bf}{\it}
{\bf}{\it}

\begin{document}
\bibliographystyle{plain}

\title{Idempotents of Hecke algebras of type $A$.}
\author{
A.K. Aiston and H.R. Morton 
	\thanks{The first author was supported by EPSRC grant GR/J72332.}\\
Department of Mathematical Sciences,\\
University of Liverpool,\\
Liverpool,\\
L69 3BX.}
\date{ \today\\Version 1.6}
\maketitle
\pagestyle{myheadings}
\markright{Idempotents of the Hecke algebra.}


\input AMsec1
\input AMsec2
\input AMsec3
\input AMsec4

\input AMsec5

\input AMsec6


\bibliography{library}
\end{document}

%% file: AMsec1
\begin{abstract}
We use a skein-theoretic version
of the Hecke algebras of type $A$ to present three-dimensional
diagrammatic views of
Gyoja's idempotent elements,
based closely on the corresponding Young diagram $\lambda$.
In this context we give straightforward calculations  for
the eigenvalues $f_\lambda$ and $m_\lambda$ of two natural 
central elements in the Hecke algebras, namely 
the full curl and the sum of the Murphy operators.
We discuss their calculation also in terms of the framing factor 
associated to the appropriate irreducible representation 
of the quantum group, $SU(N)_q$.
\end{abstract}

\section{Introduction}
The Hecke algebras of type $A$ have played a central role in the original
discovery of the Homfly polynomial for links, as described by Ocneanu 
\cite{ocneanu}.
They also relate to the representation theory of the quantum groups $SU(N)_q$ 
in a very
direct counterpart to the classical relations between the symmetric group 
algebras
and the representation theory of the linear groups $GL(N)$, \cite{wasser,wen3}.
As a result they form a natural means of transition between the invariants of 
a knot
determined from quantum-group representations and those determined by using
the Homfly polynomial of suitable decorations of the knot.

Because of the close connection with the Homfly polynomial there are 
well-established 
skein-theoretic ways of viewing the type $A$ Hecke algebras in terms of linear
combinations
of braids or tangles, modulo some simple linear relations, \cite{mt}.
Irreducible representations of the Hecke algebras are known to be 
associated to Young diagrams via quantum analogues of the symmetriser 
and anti-symmetriser.
In this paper we construct three-dimensional skein-theoretic versions
of the idempotent elements in these algebras, 
based on Gyoja's algebraic versions in \cite{gyoja}, which can be 
directly visualised  in terms of the corresponding Young diagram.

We give direct skein-theoretic proofs of some of their properties.
These lead to a simple calculation of the eigenvalues 
$f_\lambda$ and $m_\lambda$ for the full curl and for the sum of 
the Murphy operators when applied to the idempotent for the Young 
diagram
$\lambda$ in the Hecke algebra $H_n$ with $n=\vert\lambda\vert$.
While these eigenvalues are already known, \cite{dip,chak},
we believe that our approach brings out certain unexpected features.

We exploit the relationships
between the Hecke algebras and
Homfly skein theory
to express the idempotents as linear combinations of the positive permutation
braids. 

The three-dimensional idempotents described here use a natural set of 
parameters $x,v$ and $s$. 
The closure of the idempotent for the Young diagram $\lambda$
in the skein of the annulus gives a linear combination of
patterns which determines a knot invariant with parameters $x,v$ and $s$,
depending on $\lambda$. The substitution $s=\E^{h/2}$, along 
with $v=\E^{-Nh/2}=s^{-N}$ and $x=\E^{-h/2N}=s^{-1/N}$ then gives for each $N$
an invariant of the knot which is the 
$SU(N)_q$--invariant of the knot when coloured by the irreducible 
representation
of $SU(N)_q$ with the same Young diagram $\lambda$, 
where $q=\E^h$. This gives an explicit way to calculate the $SU(N)_q$--
invariants of the knot in terms of the Homfly polynomials of its satellites.

Yokota \cite{yoko} has described a rather similar skein-theoretic version of 
the idempotents, although he does this individually for each $N$ in terms of a
single parameter, while here we make use
of the extra variables $x$ and $v$ from the Homfly
skein to handle all values of $N$ at once.  

Akutsu, Deguchi and Wadati \cite{daw,dwa} have described idempotents in the
cases of individual $SU(N)_q$ which can be used to give knot 
invariants based on a chosen Young diagram, as well as a means
for extending these to 2-variable invariants, which are essentially these
Homfly based satellite invariants. 

There are accounts by Wenzl \cite{wenbcd,wen1} both in this case and 
in the setting 
of the quantum groups of types $B,C$ and $D$, of the general theory 
connecting quantum invariants 
to 2-variable knot invariants using idempotents of suitable algebras.

Our motivation is to show how the skein theory approach, based on the Homfly
polynomial, can be used very simply and explicitly in dealing with certain 
aspects
of the quantum $SU(N)$ invariants. In many cases the quantum and skein
techniques are complementary, and it is useful to have a better understanding 
of the  transitions between them, so as to translate readily from one context 
to the other.

For example, the eigenvalues $f_\lambda$ for the full curl give the
 eigenvalues for the Casimir operators of $SU(N)_q$, and thus the framing 
factors, in a way which
focuses attention on the nature of their dependence on $N$. Equally either 
version of knot invariants gives a convenient starting place for the 
Witten-Reshetikhin-Turaev $SU(N)$ manifold invariants, which can be 
constructed by suitable substitution of a root of unity in a combination
of link invariants, \cite{knot96}. The facility to change readily between 
the quantum and 
skein  invariants then allows for flexibility in interpretation and 
a consequent 
anticipation of potentially interesting features.

%% file: AMsec2
\section{The Hecke algebra and Homfly skein theory}
\label{iso}

We give a brief description of skein theory based on planar pieces
of knot diagrams and a framed version of the Homfly polynomial.  
The ideas go back to Conway and have been substantially developed
by Lickorish and others.  A fuller version of this account can be 
found in \cite{nato}.  At a later stage we shall expand our view from 
diagrams to actual pieces  of knot  lying in controlled regions of
$3$-dimensional space, under suitable equivalence.

\begin{figure}
\[
x^{-1}\,{\cal X}\left(\raisebox{-1mm}{\,\epsfxsize.2in\epsffile{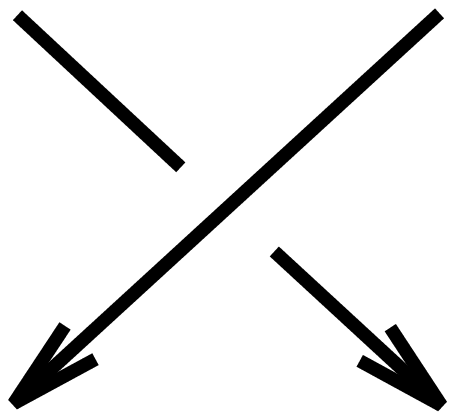}
								\,}\right)
     \,-\, x\,{\cal X}\left(\raisebox{-1mm}{\,\epsfxsize.2in\epsffile{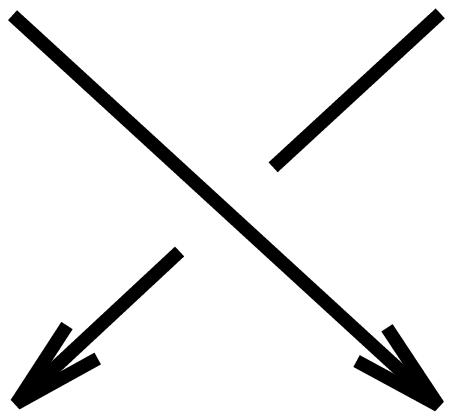}
								\,}\right)
        \,=\,z\,{\cal X}\left(\raisebox{-2mm}{\,
			\epsfxsize0.2in\epsffile{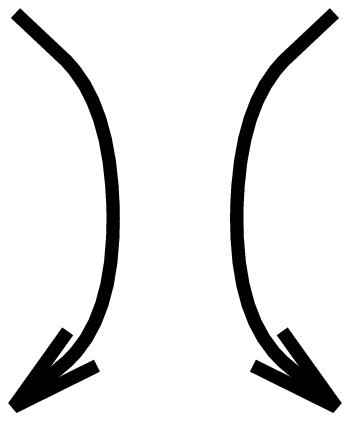}\,}
                                                        \right)
\]
\[
{\cal X}\left(\raisebox{-3mm}{\,\,\epsfxsize0.2in\epsffile{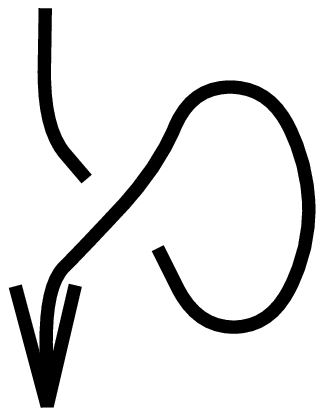}\,}\right)
        \,=\,(xv^{-1})\,{\cal X}\left(
                        \raisebox{-3mm}{\,\,\epsfysize0.35in\epsffile{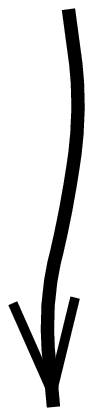}\,\,}
                                        \right)\;.
\]
\caption[]{The skein relations for the framed Homfly polynomial}
\label{homf}
\end{figure}

We shall work with the framed Homfly polynomial
${\cal X}$. This is an invariant of framed oriented links,
 constructed from the Homfly polynomial by setting 
${\cal X}(L)=(xv^{-1})^{\omega (D)}P(L)$. Here $\omega (D)$
is the writhe (the sum of the signs of the crossings) of any
 diagram $D$ of the framed link $L$ which realises the chosen framing
by means of the `blackboard parallel'. It is determined, up to a scalar,
by the skein relations in Fig.~\ref{homf}
for conventionally framed link diagrams which differ only as shown.
We will normalise ${\cal X}$ to take the value $1$ on the
empty knot, taking $\displaystyle{{\cal X}={v^{-1}-v\over z}}$ for the unknot
with zero framing.

\begin{figure}
\[
\begin{array}{ccc}
x^{-1}\quad\raisebox{-1mm}{\epsfxsize.2in\epsffile{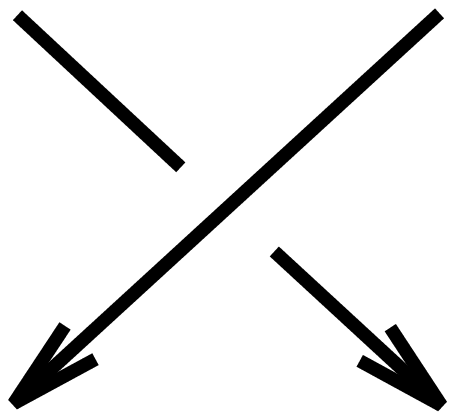}}\quad
        -\quad x\quad\raisebox{-1mm}{\epsfxsize.2in\epsffile{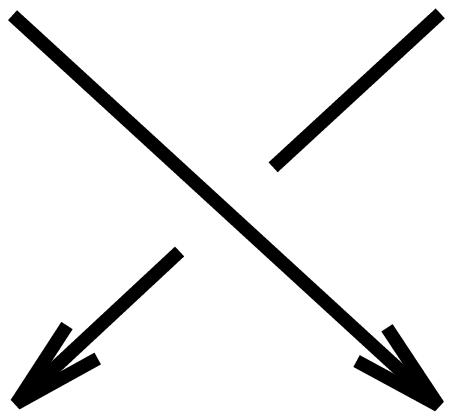}}\quad
=\quad z \quad \raisebox{-2mm}{\epsfxsize.2in\epsffile{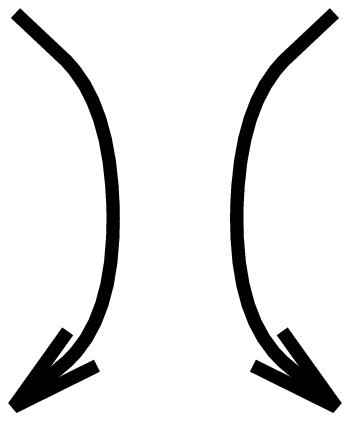}}\;,
&\quad
\mbox{\phantom{and}}\quad
&
\raisebox{-3mm}{\,\,\epsfxsize0.2in\epsffile{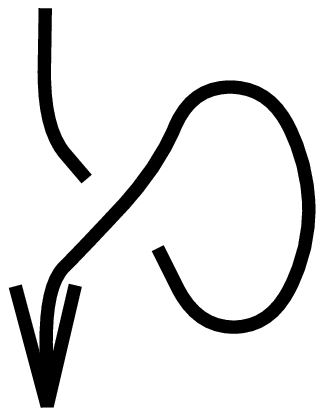}\,}
        \quad=\quad (xv^{-1})\,\,
                \raisebox{-3mm}{\,\,\epsfysize.3in\epsffile{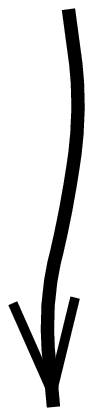}\,\,}\quad.
\end{array}
\]
\caption{The skein relations for the Homfly skein of a surface}
\label{skein}
\end{figure}

\begin{figure}
 \[
D\,\sqcup\,\raisebox{-2mm}{\epsfysize.3in\epsffile{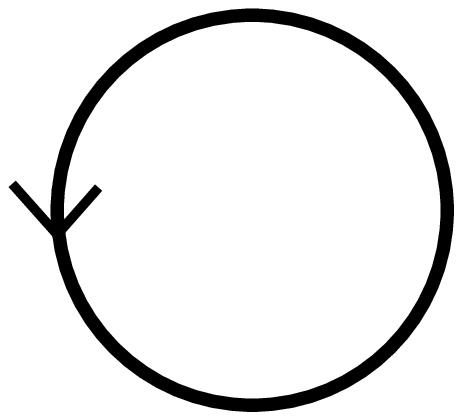}}\qquad
                                        =\qquad {v^{-1}-v\over z}\quad D\;,
\]
\caption{The skein relation for a null homotopic loop}
\label{null}
\end{figure}

Let $F$ be a planar surface and fix  a (possibly
empty) set of distinguished points on the boundary.   
We consider diagrams in $F$, consisting of oriented arcs joining any
distinguished boundary points and oriented closed curves, up to Reidemeister moves
II and III. They carry the
implicit framing defined by the parallel curves in the diagram. Define the
{\it framed Homf\,ly skein\/} of $F$, denoted by ${\cal S}(F)$, to consist of
linear combinations of diagrams in $F$ modulo
the skein relations in Fig.~\ref{skein}
As a consequence, the relation in Fig.~\ref{null} holds in ${\cal S}(F)$,
for any diagram $D$, where the oriented circle denotes a null homotopic loop. 

We are interested in  three specific cases, namely when $F$ is the whole  plane
${\Bbb R}^2$, the annulus $S^1\times I$
or the rectangle $R^n_n \cong I\times I$ with $n$ distinguished  points
on its top and bottom edge.
In the last case we insist
that any arcs in $R^n_n$ enter at the top and leave at the bottom.
Diagrams in $R^n_n$ are termed {\it oriented\/} $n$-{\it tangles\/}, and include the case of
$n$-string braids.

A {\it positive permutation braid\/} is defined for each permutation
$\pi\in S_n$.
It is the $n$-string braid, $\omega_\pi$,
uniquely determined by the properties 
\newcounter{anna}
\begin{list}{\roman{anna})}{\usecounter{anna}}
\setlength{\rightmargin}{\leftmargin}
\item all strings are oriented from top to bottom
\item for $i=1,\ldots,n$ the $i$th string joins the point numbered $i$ at the 
top of the braid to the point numbered $\pi(i)$ at the bottom of the braid,
\item all the crossings occur with positive sign and each pair of strings cross at most once.
\end{list}
These were first defined by
Elrifai and Morton \cite{fai}.
We can think of the braid strings as sitting in layers, with the 
first string at the back and the $n$th string at the front.

We define the {\it negative permutation braid\/} for $\pi$ 
in exactly the same manner as the positive permutation braid except that
we demand that all the crossing be negative.  We shall denote
this braid by $\overline{\omega}_\pi$.
The inverse of $\omega_\pi$ is the negative permutation braid
with permutation $\pi^{-1}$, thus
$\omega_\pi^{-1}=\overline{\omega}_{\pi^{-1}}$.

The skein  ${\cal S}({\Bbb R}^2)$ is just the set of
linear combinations of
framed link diagrams, modulo the skein relations. Every diagram $D$ represents a scalar
multiple of the empty diagram, with the multiple being simply  ${\cal X}(D)$.

It is shown in \cite{mt} that the skein ${\cal S}(R_n^n)$ is spanned
by the $n!$ positive permutation braids, and that these are linearly independent.
The elementary braid $\sigma_i$, which is the positive permutation braid
for the transposition $(i\,i+1)$, satisfies the relation $x^{-1}\sigma_i-x
\sigma_i^{-1}=z$ in the skein. The skein forms an algebra over $\Lambda$ with 
multiplication derived from the concatenation of diagrams. As is conventional for braids, we
write $ST$ for the diagram given by placing  diagram $S$ above diagram $T$.
The resulting algebra is a quotient of the braid-group algebra. It is shown in \cite{mt}
that the algebra ${\cal S}(R_n^n)$ is isomorphic to the Hecke algebra $H_n$ of type $A$,
with the explicit presentation
\[
H_n\quad=\quad
\left<
\begin{array}{ccc}
        \begin{array}{ccc}
        \sigma_i& : & i=1,\ldots, n-1\\
               &   &
        \end{array}
&\left.\begin{array}{c}
         \\
         \\
        \end{array}     \right\vert

&       \begin{array}{l}
        \sigma_i\sigma_j=\sigma_j\sigma_i~:~\vert i-j\vert>1\\
        \sigma_i\sigma_{i+1}\sigma_i=\sigma_{i+1}\sigma_i\sigma_{i+1}\\
	x^{-1}\sigma_i-x\sigma^{-1}_i=z\;,
        \end{array}
\end{array}
\right>\;.
\]

There are various presentations of the Hecke algebra in the literature.
Here we have used a coefficient ring $\Lambda$ with 3 variables $x,v$ and
$z$.  The variable $v$ is needed in the skein
when we want to write a general tangle in terms of the basis of permutation braids,
but it does not appear in the relations. The variable $x$ keeps track of the writhe
of a diagram, and can be dropped without affecting the algebraic properties. Define
an algebra  $H_n(1,z)$,  obtained
from $H_n$ by setting $x=1$, in terms of generators $\rho_i$ where the 
quadratic relation is given by $\rho_i-\rho_i^{-1}=z$.

In algebra texts such as \cite{chak,dip} 
the Hecke algebra is usually presented with a quadratic
relation whose roots are $q$ and $-1$. We will denote this presentation by $H_q(n)$. 
With generators  $\tau_i$, the quadratic relation is
$\tau_i^2=(q-1)\tau_i+q$.
Note that these presentations are all isomorphic.
\[
\begin{array}{ccccc}
H_n &\cong&H_n(1,z)&\cong&H_q(n)\\
\sigma_i&\mapsto&x\rho_i& &\\
	&	&\rho_i &\mapsto&s^{-1}\tau_i\\
\end{array}
\]
where we relate $q$ and $z$ by setting $z=s-s^{-1}$ and $q=s^2$.

In what follows we shall generally use the skein form ${\cal S}(R_n^n)$ of
$H_n$, but the results can be transferred immediately 
to the other versions by these isomorphisms.
\medskip

A {\it wiring\/} $W$ of a surface $F$ into another surface $F'$
is a choice of inclusion of $F$ into
$F'$ and a choice of a fixed diagram of curves and arcs in $F'- F$ whose
 boundary  is the union of the distinguished
sets of $F$ and $ F'$. A wiring $W$ determines naturally a $\Lambda$-linear map 
${\cal S}(W):{\cal S}(F)\to {\cal S}(F')$.

\begin{figure}
\[
\epsfysize.6in\epsffile{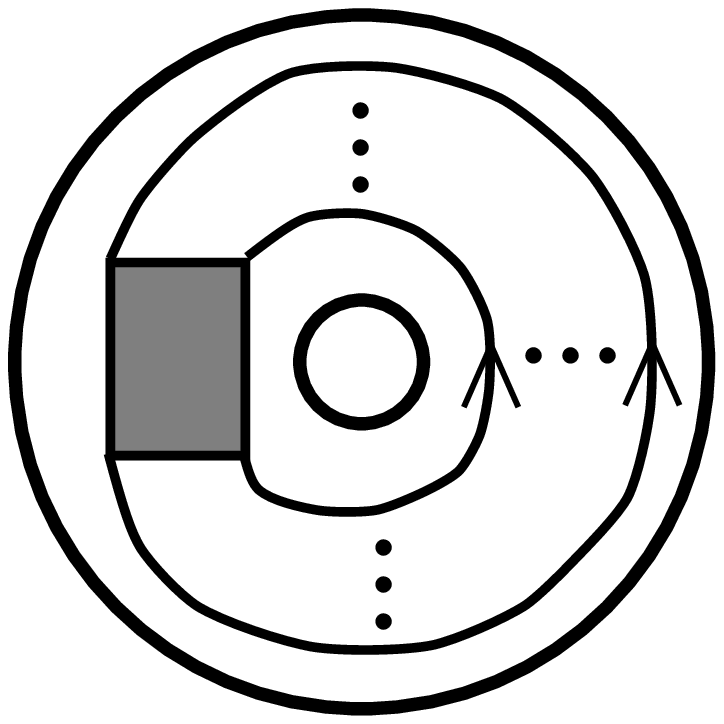}
\]
\caption[]{The wiring of ${\cal S}(R_n^n)$ into ${\cal S}(S^1\times I)$}
\label{wire}
\end{figure}

We can wire the rectangle $R^n_n$
into the annulus as indicated in Fig.~\ref{wire}.
The resulting diagram in the annulus is called the {\it closure\/} of the oriented
tangle.
We shall also use the term `closure' for the  $\Lambda$-linear map
from each Hecke algebra $H_n$ to the skein of the annulus induced by this wiring.

The skein ${\cal S}(S^1\times I)$ of the annulus itself forms an algebra, whose
 product is given by
stacking the annuli one inside the other.  This product is
obviously commutative (lift the inner annulus up and stretch it so
that the outer one will fit on the inside of it).
Write ${\cal C}$ for ${\cal S}(S^1\times I)$ regarded as a  $\Lambda$-algebra 
in this way.
Let ${\cal C}^+$ be the sub-algebra spanned by the closures of oriented tangles.
\label{bazza}
Turaev \cite{turbas} showed that ${\cal C}^+$ is freely generated
as an algebra
by $\{ A_m$, $m\in {\Bbb N}\}$,  where $ A_m$ is the closure of
the positive permutation braid for the cycle $(1\,2\,\ldots\,m)$.

Write ${\cal C}^{(n)}$ for the linear subspace spanned by the closures
of oriented $n$-tangles.
This is spanned by all terms of the form
$( A_{i_1})^{j_1}( A_{i_2})^{j_2}\cdots( A_{i_p})^{j_p}$
where $i_k, j_k\in{\Bbb N}$ and $\sum_{k=1}^pi_kj_k=n$. Hence the algebra
${\cal C}^+$ is graded as
\[  
{\cal C}^+=\bigoplus_{n=0}^\infty {\cal C}^{(n)}\;,
\]
and each subspace ${\cal C}^{(n)}$ is the image of the Hecke algebra $H_n$
under the closure map.

%% file: AMsec3
\section{Young diagrams}
\label{bib}
There is a wealth of detail about the features of Young
tableaux in many texts such as \cite{weyl,repthry,jones}.
Here we emphasize certain properties which will be to the
fore in this article.

A partition of $n$ can be represented by a {\it Young diagram\/};
a collection of $n$ cells arranged in
rows, with $\lambda_1$ cells in the first row, $\lambda_2$ cells 
in the second row up to $\lambda_k$ cells in the $k$th row where
$\lambda_1\geq\lambda_2\geq\cdots\geq\lambda_k>0$
and $\sum_{i=1}^k\lambda_i=n$.
We shall denote both the partition and
its Young diagram by $\lambda$. The Young diagram
for $(0)$ is the empty diagram.
We denote the number of cells
in $\lambda$ by $\vert\lambda\vert$.
The {\it conjugate\/} $\lambda^\vee$ of $\lambda$
is the Young diagram whose rows form the columns of $\lambda$.
Any cell for which a legitimate Young diagram remains
after it has been removed will be called an {\it extreme cell\/}.
To each extreme cell we associate an {\it extreme rectangle\/},
namely those cells above and to the left of it in the Young diagram.
When there is a cell in the 
$i$th row and $j$th column of $\lambda$ we write $(i,j)\in \lambda$, and refer
to $(i,j)$ as the {\it coordinates\/} of the cell.

\begin{figure}
\[
\nu\quad=\quad\raisebox{-3mm}{\epsfxsize.4in\epsffile{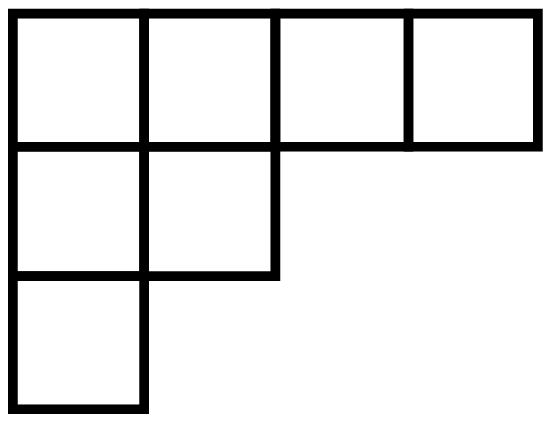}}
\qquad\qquad
\nu^\vee\quad=\quad\raisebox{-5mm}{\epsfxsize.3in\epsffile{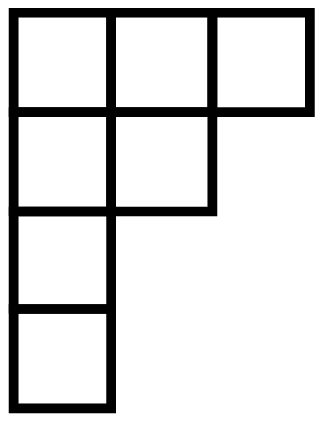}}
\]
\caption{The Young diagram $\nu=(4,2,1)$ and its conjugate}
\label{nu}
\end{figure}

\begin{figure}
\[
\epsfxsize.4in\epsffile{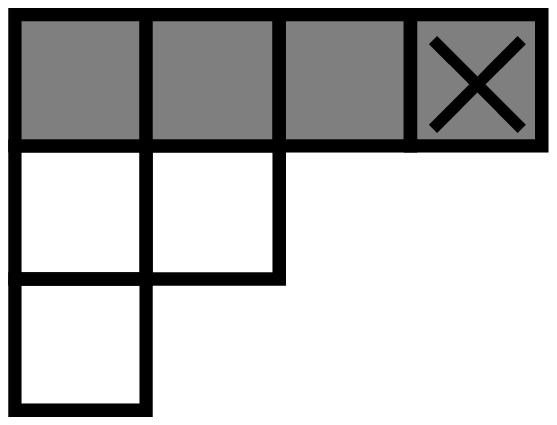}\qquad
\epsfxsize.4in\epsffile{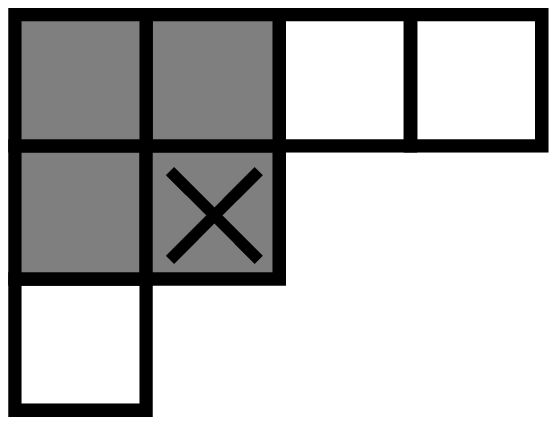}\qquad
\epsfxsize.4in\epsffile{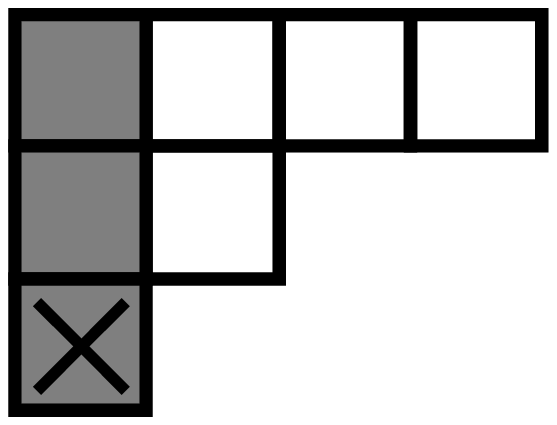}
\]
\caption[]{The extreme cells and rectangles of $\nu$}
\label{treme}
\end{figure}

We will work with the example $\nu=(4,2,1)$ throughout
this paper. 
The conjugate of $\nu$ is $\nu^\vee=(3,2,1,1)$ and 
$\vert\nu\vert=\vert\nu^\vee\vert=7$.
The Young diagrams for $\nu$ and  $\nu^\vee$ are shown in Fig.~\ref{nu}.
There are three extreme cells in $\nu$, marked by a cross in Fig.~\ref{treme}
with their associated extreme rectangles shaded in.  
The coordinates of the extreme cells are $(1,4)$, $(2,2)$ and $(3,1)$ 
respectively.

Let $T'(\lambda)$ be an
assignment of the numbers $1$ to $n$ to the cells of $\lambda$ such that
the numbers increase from left to right along the rows and top to bottom
down the columns. We call $T'(\lambda)$ a {\it standard tableau\/}.
In particular $T(\lambda)$ will denote the tableau where the cells of the
Young diagram are numbered from $1$ to $n$ along the rows.

\begin{figure}
\[
T(\nu)\ =\ \raisebox{-3mm}{\epsfxsize.533in\epsffile{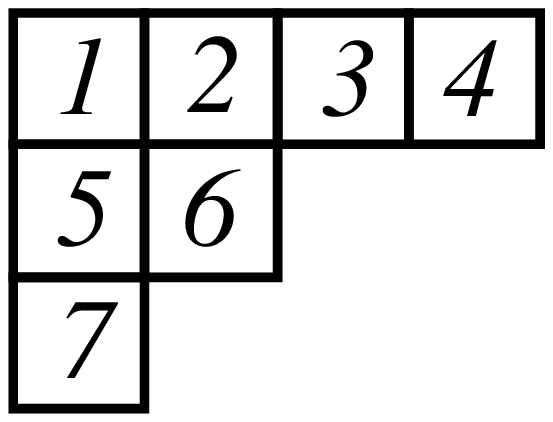}}
\qquad\qquad
T(\nu^\vee)\ =\ \raisebox{-5mm}{\epsfxsize.4in\epsffile{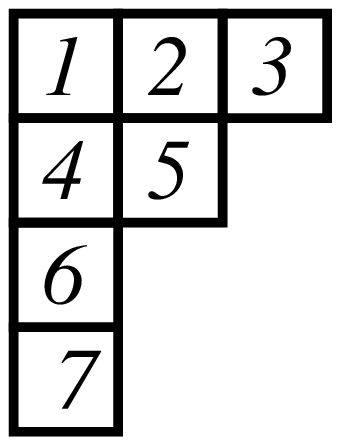}}
\]
\caption[]{The tableaux $T(\nu)$ and $T(\nu^\vee)$}
\label{eaux}
\end{figure}

Note that the transposition of rows and columns doesn't take $T(\lambda)$ 
to $T(\lambda^\vee)$. We define the permutation 
$\pi_\lambda$ by $\pi_\lambda(i)=j$
where the transposition of $\lambda$
carries the cell $i$ in $T(\lambda)$
to the cell $j$ in $T(\lambda^\vee)$. 
The tableaux $T(\nu)$ and $T(\nu^\vee)$ are shown in Fig.~\ref{eaux}
and $\pi_\nu=\left(2\,4\,7\,3\,6\,5\right)$.

We now give some standard combinatorial results 
about permutations and Young diagrams, in a context which extends readily
to positive permutation braids and the Hecke algebras.

Let $\lambda$ and $\mu$ be  Young diagrams with 
$\vert\lambda\vert=\vert\mu\vert=n$.
We say that $\pi\in S_n$ {\it separates\/} $\lambda$ from $\mu$ if no pair
of numbers in the same row of $T(\lambda)$ are mapped by $\pi$ to 
the same row of $T(\mu)$. 
The permutation $\pi_\lambda$, for example, separates $\lambda$
from its conjugate $\lambda^\vee$.

Write $R(\lambda)\subset S_n$ for the subgroup of permutations which 
preserve the rows of $T(\lambda)$. Each $R(\lambda)$ is generated by some
subset of the elementary transpositions $(i\,i+1)$. For example  $R(\nu)$ is
 generated by 
$\{(12),(23),(34),(56)\}$. 

It is easy to see that if $\pi$ separates $\lambda$ from $\mu$ then so does
$\rho\pi\sigma$ for any $\rho\in R(\lambda),\sigma\in R(\mu)$. Conversely, 
it can be shown that if $\pi$ separates $\lambda$ from $\lambda^\vee$ then 
$\pi=\rho\pi_\lambda\sigma$
with $\rho\in R(\lambda)$ and $\sigma\in R(\lambda^\vee)$.
We say that $\lambda$ is {\it just separable\/} from      
$\lambda^\vee$. If no permutation $\pi\in S_n$
separates $\lambda$ from $\mu$ then we
call $\lambda$ and $\mu$ {\it inseparable\/}.
\medskip

Order the Young diagrams by lexicographical ordering of their rows. Thus
 $\lambda>\mu$ when there exists $t$ with  $\lambda_t>\mu_t$ 
and $\lambda_i=\mu_i$ for $i<t$.

\begin{aslemma}
\label{sep}
If $\lambda>\mu$ then $\lambda$ and $\mu^\vee$ are inseparable.
\end{aslemma}
\begin{proof}
By induction on $t$.
Suppose first that $t=1$, and let $\pi\in S_n$. 
Since $\mu^\vee$ has $\mu_1<\lambda_1$ rows, $\pi$ must map at least two
numbers from the first row of $\lambda$ to the same row of $\mu^\vee$. Hence
$\pi$ does not separate $\lambda$ from $\mu^\vee$.

When $t>1$,  consider the Young diagrams $r(\lambda)$ and $r(\mu)$,
where $r(\alpha)$ is the Young diagram
obtained from $\alpha$ by removing the first row.
We have $r(\lambda)>r(\mu)$, and the diagrams first differ in the $(t-1)$st row.
Hence, by the induction hypothesis, $r(\lambda)$ and  $r(\mu)^\vee$ are inseparable.
It remains to prove that if $\lambda$ were separable from $\mu^\vee$ 
then $r(\lambda)$ would be separable from $r(\mu)^\vee$, 
giving us a contradiction.

Suppose that $\pi$  separates $\lambda$ from $\mu^\vee$. Since
$\lambda_1=\mu_1$, the number of cells in the first row of $\lambda$ is equal
to the number of rows of $\mu^\vee$. Then $\pi$ must send exactly one
cell from the first row of $\lambda$ to each row of $\mu^\vee$.  We can suppose,
without loss of generality, that it is the first cell in each row of $\mu^\vee$,
for if not, there is a transposition $\tau\in R(\mu^\vee)$ 
which will switch the first cell of
the row with the image of the cell in the first row of $\lambda$, and $\pi \tau$
will also separate $\lambda$ from $\mu^\vee$. 
Restrict $\pi$ to all but the first row of $\lambda$. Its image
will be exactly the cells of $r(\mu)^\vee$, and the restricted permutation
 separates $r(\lambda)$ from $r(\mu)^\vee$.
\end{proof}
\begin{ascor}
If $\lambda\neq\mu$ then either
$\lambda$ and $\mu^\vee$  are inseparable  or 
$\lambda^\vee$ and $\mu$ are inseparable.
\end{ascor}
\begin{aslemma}
\label{pimult}
Let $\pi$ be a permutation which separates $\lambda$ from $\lambda^\vee$. 
Then the positive permutation braid $\omega_\pi$ can be written as
$\omega_\pi=\omega_{\rho}\omega_{\pi_\lambda}\omega_{\sigma}$ for some 
$\rho\in R(\lambda)$ and $\sigma\in R(\lambda^\vee)$.
\end{aslemma}
\begin{proof}
We already know that we can write $\pi=\rho\pi_\lambda\sigma$. It is enough
to show that $\omega_{\rho}\omega_{\pi_\lambda}\omega_{\sigma}$ 
is a positive permutation braid, for it is then determined 
as $\omega_\pi$ by its permutation $\pi$.

The key feature of $\omega_{\pi_\lambda}$ is that pairs of strings which start
in the same row of $T(\lambda)$ or finish in the same row of $T(\lambda^\vee)$ 
don't cross.
The only pairs which cross in $\omega_\rho$ or in $\omega_\sigma$ are in the 
same row of $T(\lambda)$ or $T(\lambda^\vee)$ respectively. Consequently each
pair of strings in the product braid $\omega_{\rho}\omega_{\pi_\lambda}
\omega_{\sigma}$ crosses at most once. All  crossings are positive, so this is 
a positive permutation braid.
\end{proof}

This is a special case of a property of positive permutation braids which 
will be useful later. 
\begin{aslemma}
In each right coset $R(\lambda)\pi'$ of the subgroup
$R(\lambda)$ there is a unique $\overline{\pi}$ such that 
$\omega_\pi=\omega_\rho \omega_{\overline{\pi}}$,
where $\rho\in R(\lambda)$ and $\pi=\rho\overline{\pi}$ runs through the coset.
\end{aslemma}
\begin{proof} 
Choose $\overline{\pi}$ so that no two strings which start in the same
row of $\lambda$ cross. Then $\omega_\rho \omega_{\overline{\pi}}$ is a
 positive permutation braid.
\end{proof}
\begin{ascor} 
\label{picoset}
If $\pi=\rho\pi'$ with $\rho\in R(\lambda)$ then
$\omega_\pi=\omega_{\rho'}\overline{\omega}_{\rho''} \omega_{\pi'}$, for some
$\rho',\rho''\in R(\lambda)$.
A similar result holds for left cosets.
\end{ascor}

%% file: AMsec4
\section{Idempotents}
\label{idemp}

 The Hecke algebra $H_n$ is closely related to the
group algebra of $S_n$, whose idempotents  are described by the classical Young
symmetrisers. For a  Young diagram $\lambda$ its  Young symmetriser is the
product of the sum of permutations which preserve the rows of $T(\lambda)$ and
the alternating sum of permutations which preserve columns. With care it is
possible to make a similar construction of idempotents
 in $H_n$, replacing permutations by suitably weighted positive permutation
braids. Jones \cite{hecke} gives a good description of the two idempotents
corresponding to single row and  column Young diagrams. Other authors, for
example Wenzl and Cherednik, have given descriptions for general $\lambda$, but
we shall here adapt the construction of Gyoja \cite{gyoja} to construct
idempotents in $H_n$ regarded as the skein ${\cal S}(R^n_n)$.  We shall follow 
the account  in \cite{nato} for the basic row and column idempotents, and use
these to  construct an idempotent for each Young diagram $\lambda$.

We start from the visually appealing $3$-dimensional picture for the idempotent 
as a linear combination of braids in a $3$-ball based very closely on the
diagram $\lambda$  rather than as a linear combinations of diagram  in a
rectangle. 

We consider a $3$-ball $B\cong B^3$, with a chosen subset $P$
of $2n$ points on its boundary sphere, designated as $n$ inputs $P_I$
and $n$ outputs $P_O$.  An
oriented tangle $T$ in $(B,P)$ is made up of $n$ oriented arcs in $B$ joining
the points $P_I$ to the points $P_O$, together with any number of oriented
closed curves. The arcs and curves of $T$ are assumed to carry a framing
defined by a specific choice of parallel for each component.

The skein ${\cal S}(B,P)$ is defined as linear combinations of such tangles,
modulo the framed Homfly skein relations applied to tangles which differ only
as in Fig.~\ref{skein} inside some ball. The case when $B=D^2\times I$ and 
the points $P_I$ and $P_O$
are lined up along the top and bottom respectively,
gives a skein which can readily be identified with ${\cal S}(R_n^n)=H_n$. 
There is a homeomorphism mapping any other pair
$(B',P')$ to this pair, when $|P'|=2n$. 
This induces a linear isomorphism from each ${\cal S}(B',P')$ to 
 the Hecke algebra $H_n$.

As in the case of diagrams, a {\it wiring} $W$ of $(B,P)$ into
$(B',P')$ is an inclusion of the ball $B$ into the interior of $B'$ and a
choice of framed oriented arcs in $B-B'$ ending with compatible orientation at
the boundary points $P\cup P'$.  Given a tangle $T$ in $(B,P)$ and a
wiring $W$ their union determines a tangle $W(T)$ in $(B',P')$, 
and induces a linear map ${\cal S}(W):{\cal S}(B,P) \to {\cal S}(B',P')$.

The region between two balls is homeomorphic to $S^2\times I$. A simple example
of wiring $W$ consists of $n$ arcs with each lying monotonically in the $I$
coordinate, sometimes called an $n$-braid in $S^2$. In such a case the map
${\cal S}(W)$ is a linear isomorphism whose inverse is induced by the inverse
braid. Such a wiring can always be chosen to determine an explicit isomorphism
from any ${\cal S}(B,P)$ to $H_n={\cal S}(R_n^n)$ when $|P|=2n$.

In terms of calculation the $3$-dimensional viewpoint has advantages and
disadvantages. There  may not be a natural way to compose tangles in the skein
${\cal S}(B,P)$,
so the isomorphism with $H_n$ does not immediately carry any
algebra information. On the other hand  a suitable choice of wiring may show 
when certain elements are zero in the skein, and hence exhibit relations 
in $H_n$.

The following picture gives the heart of the construction of the idempotent for
the Young diagram $\lambda$. It lies in the skein of $B^3\cong D^2\times I$
where the points $P_I$ and $P_O$  are the centres of cells of    templates  in the shape
 of $\lambda$ at the top and
bottom respectively.  The strings in each rows are first grouped
together using a linear combination $a_j$ of braids for a row with $j$ cells.
Below this the strings in  the columns are grouped with linear combinations
$b_j$ of braids, to define an element $E_\lambda$ of the skein. The elements
$a_j$ and $b_j$ are Jones'  basic row and column quasi-idempotents, described
shortly  in more detail. 

\begin{figure}
\[
\epsfxsize1.75in\epsffile{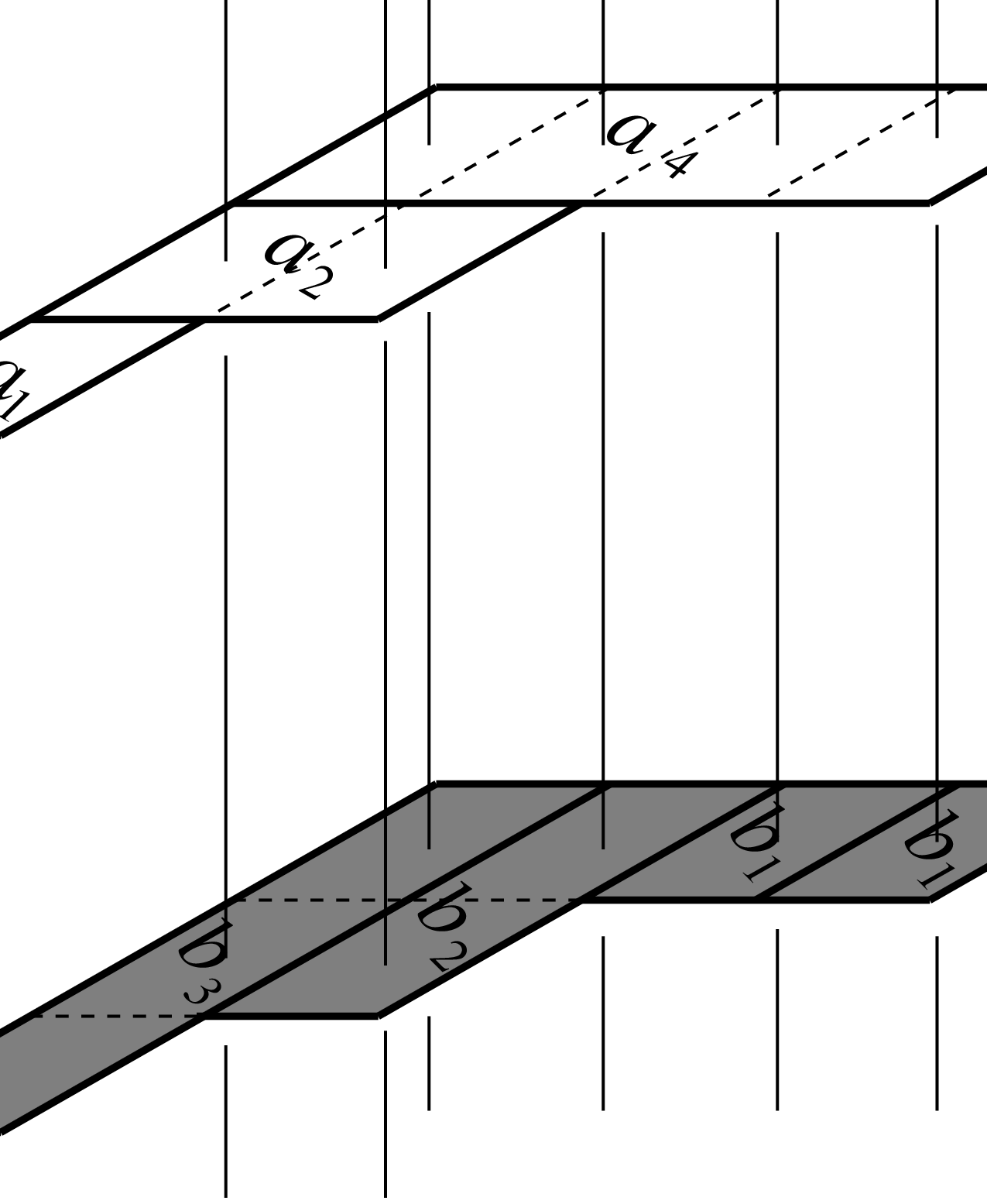}
\]
\caption[]{The $3$-dimensional quasi-idempotent associated to $\nu$}
\label{3D}
\end{figure}

We give the $3$-dimensional picture for $E_\nu$ in Fig.~\ref{3D}.  Note how the
strings are first grouped in rows and then in columns.  Although the braids run
vertically the boxes have been drawn in the horizontal plane to emphasize how
the idempotent is related to the shape of the Young diagram. 

\begin{figure}
\[
\epsfxsize1in\epsffile{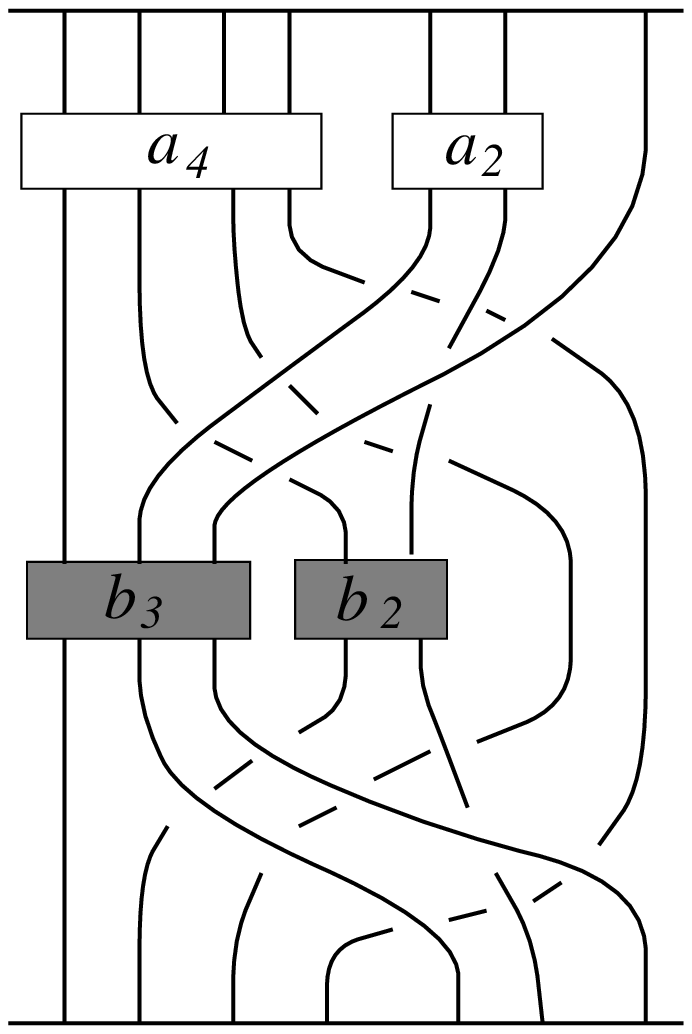}
\]
\caption[]{The quasi-idempotent $e_\nu$}
\label{hats}
\end{figure}

This $3$-dimensional view will be discussed further in the next section, along
with some closely related choices of wiring.   For the moment we
will work with a
$2$-dimensional picture given by using the standard tableau $T(\lambda)$ to
produce a wiring which  straightens  out the strings to lie in order on a line.
The corresponding figure, more immediately related to Gyoja's idempotents,  is
shown in  Fig.~\ref{hats}.
\medskip

We first define the row and column elements $a_j$ and $b_j$, following
the account  in Morton \cite{nato}.

Write $\displaystyle{
E_n(\sigma_1,\sigma_2,\cdots,\sigma_{n-1}) =\sum_{\pi\in S_n} \omega_\pi}$
for the sum of the  positive permutation braids.

\begin{astheorem}[\protect\cite{nato}]
\label{factor}
For each $i$, we can factorise $E_n$ in $H_n$ as
\[
E_n=(\sigma_i+1)E_n^{(i)}, \quad E_n={E'}_n^{(i)}(\sigma_i+1),
\]
where $\displaystyle{E_n^{(i)}=\sum_{\pi(i)<\pi(i+1)}\omega_\pi}$ and $\displaystyle{{E'}_n^{(i)}=\sum_{\pi(i)<\pi(i+1)}\omega_{\pi^{-1}}}$ .
\end{astheorem}
\begin{proof}
The proof follows the argument in \cite{nato}, adapted for the fact that our convention for labelling permutations is the opposite of that in \cite{nato}.

Given  $i$, we can pair the 
permutations as follows. 
For each permutation $\pi $ consider its 
composite $\pi '=\pi \circ (i\  i{+}1)$ with the 
transposition $(i\  i{+}1)$. Exactly one of the 
pair preserves the order of $i$ and $i+1$. 
Suppose that it is $\pi $, so that  $\pi (i)<\pi 
(i+1)$ and $\pi'(i)>\pi'(i+1)$. Then the braid $\sigma _{i}\omega_{\pi }$ 
 is itself a 
positive permutation braid. Since its permutation 
is $\pi '$ we have $\sigma _{i}\omega_{\pi }=\omega_{\pi 
'}$. Then

\begin{eqnarray*}E_{n}&=&\displaystyle \sum_{\pi (i)<\pi (i+1)} 
\omega_{\pi } \quad+ \displaystyle \sum_{\pi '(i)>\pi 
'(i+1)}\omega_{\pi '} \\
&=&\displaystyle \sum_{\pi (i)<\pi (i+1)} \omega_{\pi } 
\quad+ \displaystyle \sum_{\pi (i)<\pi (i+1)}\sigma _{i}\omega_{\pi 
} \\ 
&=&(\sigma _{i}+1)E_{n}^{(i)},
\end{eqnarray*}
with $\displaystyle E_{n}^{(i)}=\sum_{\pi 
(i)<\pi (i+1)}\omega_{\pi } $.

The second equation in the  theorem follows from the first by applying the  braid-reversing antiautomorphism ${\rm Rev}$ of $H_n$, defined by 
\begin{eqnarray*}{\rm Rev}(AB)&=&{\rm Rev}(B){\rm Rev}(A)\\ 
{\rm Rev}(\sigma_i)&=&\sigma_i.
\end{eqnarray*}
Then ${\rm Rev}(\omega_\pi)=\omega_{\pi^{-1}}$, and so  \[E_n={\rm Rev}(E_n)={\rm Rev}(E_n^{(i)})(\sigma_i+1)={E'}_n^{(i)}(\sigma_i+1).\]

\end{proof}
The elementary braids $\sigma_i$ satisfy the
 quadratic relation $
(\sigma_i-a)(\sigma_i-b)=0\,$ in $H_n$,
where $a=-xs^{-1}$ and $b=xs$, with the variants $-s^{-1},s$ or $-1,q$ for
the roots $a$ and $b$ in the
other versions of the Hecke algebra.  
Define $a_n$ and $b_n$ by substituting $-a^{-1}\sigma_i$ or $-b^{-1}\sigma_i$
respectively for $\sigma_i$ in $E_n$. Thus 
\[
a_n=\sum_\pi (-a)^{-l(\pi)}\omega_\pi\;,
\quad b_n=\sum_\pi (-b)^{-l(\pi)}\omega_\pi
\]
where $l(\pi)$ is the writhe of $\omega_\pi$, known in algebraic terms as the 
length of the permutation $\pi$.
\label{lock}

\begin{ascor}
\label{scalar}
 We can factorise $a_n$ and $b_n$ in $H_n$ as
\[
a_n=(\sigma_i-a){a}_n^{(i)}={a'}_n^{(i)}(\sigma_i-a)\;, \quad
b_n=(\sigma_i-b){b}_n^{(i)}={b'}_n^{(i)}(\sigma_i-b)\;.
\] 
\end{ascor}
\begin{proof} 
Substitute $-a^{-1}\sigma_i$ or $-b^{-1}\sigma_i$
for $\sigma _i$ in Theorem~\ref{factor}.
\end{proof}

In our diagrams the elements $a_n$ and $b_n$ are drawn as rectangles.  
Those denoting  $b_n$ are shaded to distinguish them from those carrying
 $a_n$. In the classical case, with $x=s=1$, $a_n$ reduces to a sum of 
permutations and $b_n$ to an alternating sum.

Applying Cor.~\ref{scalar} and the factorisation of the
quadratic relation we obtain the following result.

\begin{aslemma}[\rm\protect\cite{nato}]
\label{linhom}
Let $\phi_a$ and $\phi_b$ be the linear homomorphisms from the Hecke algebra,
$H_n$, to the ring of scalars $\Lambda$ defined by $\phi_a(\sigma_i)=a$ and
$\phi_b(\sigma_i)=b$ for $i=1,\ldots, n-1$.
Then for all $h\in H_n$,
\[
a_nh=ha_n=\phi_b(h)a_n
\ \mbox{ and }\ 
b_nh=hb_n=\phi_a(h)b_n\;.
\]
\end{aslemma}

In particular, we note the following consequence of Lemma~\ref{linhom}.
A copy of an $a_i$ can be swallowed (from above or below) by an $a_k$, 
if $i\leq k$, at the expense of multiplying the resulting diagram by a
scalar.  This scalar is $\alpha_{(i)}$ (see page \pageref{alpha}),
which is non-zero.  Thus an $a_k$ can also throw out extra copies
of $a_i$, multiplying the resulting diagram by $\alpha^{-1}_{(i)}$.
Further, there is no net effect if we introduce 
and then later remove an $a_i$.  
This works equally well with $b$ in place of $a$.
We make further use of this property, with variants, in the next section.

\begin{figure}
\[
E_\nu(a)\quad=\quad\raisebox{-.5cm}{\epsfxsize1in\epsffile{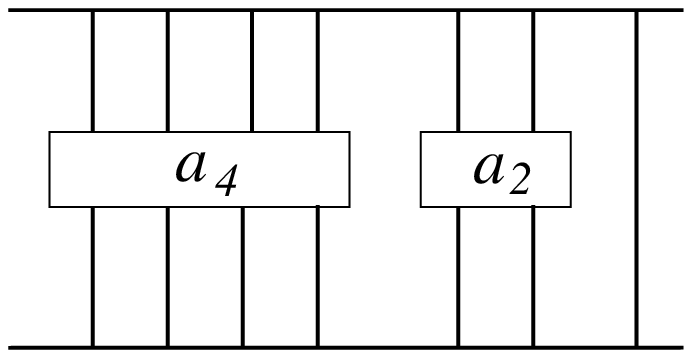}}\ .
\]
\caption{The element $E_\nu(a)$}
\label{enua}
\end{figure}

We now define the quasi-idempotent elements $e_\lambda\in H_n$ for each
Young diagram
$\lambda=(\lambda_1,\lambda_2,\ldots,\lambda_k)$.
To each cell of $\lambda$ we assign a braid string, ordered according to $T(\lambda)$.
Define $E_\lambda(a)\in H_n$ as a linear combination of braids by placing
$a_{\lambda_i}$ on the strings corresponding to the $i$th row of $\lambda$
for each $i$, and similarly $E_\lambda(b)$ using $b_{\lambda_i}$. 
The element $E_\nu(a)$ is depicted in Fig.~\ref{enua}.

Now define 
$e_\lambda
=E_\lambda(a)\omega_{\pi_\lambda} E_{\lambda^\vee}(b)\omega_{\pi_\lambda}^{-1}\in H_n$.
The element $e_\nu$ is shown in Fig.~\ref{hats}.
This picture can be obtained from the $3$-dimensional
picture of Fig.~\ref{3D} by sliding the rows apart at the top of the diagram 
and sliding the columns apart at the bottom, using the standard tableau $T(\nu)$
to determine how to order the strings.

For each $\lambda$ write $H(\lambda)\subset H_n$ for the subspace spanned by
$\{\omega_\rho, \rho\in R(\lambda)\}$. Then $E_\lambda(a),E_\lambda(b)\in H(\lambda)$
and $E_\lambda(a)h=hE_\lambda(a)=\phi_b(h)E_\lambda(a)$ for $h\in H(\lambda)$;
a similar result holds with $a$ and $b$ interchanged.
\begin{aslemma}
\label{=zero}
Let
$\lambda$ and $\mu$ be Young diagrams  with $n$ cells, and let $\pi\in S_n$
be a permutation which does not separate $\lambda$ from $\mu$.
Then
\[
E_\lambda(a) \omega_\pi E_\mu(b)=0=E_\lambda(b) \omega_\pi E_\mu(a).
\]
\end{aslemma}
\begin{proof}
Since $\pi$ does not separate $\lambda$  from $\mu$, there are two cells in some row
of $\lambda$, the $l$th say,  which are sent to two cells in the same row of 
$\mu$, the $p$th say, by $\pi$.
We can  find $\rho\in R(\lambda)$
and $\sigma\in R(\mu)$ such that $\pi'=\rho\pi\sigma$
 sends the two adjacent cells $i$ and $i+1$ in the $l$th row of $\lambda$
to cells $j$ and $j+1$ in the $p$th row of $\mu$.
By Cor.~\ref{picoset} we can write 
$\omega_\pi=\omega_{\rho'} \overline{\omega}_{\rho''}
\omega_{\pi'}\overline{\omega}_{\sigma'}\omega_{\sigma''}$, with $\rho',\rho''\in R(\lambda)$
and $\sigma',\sigma''\in R(\mu)$. Hence $\omega_\pi=h\omega_{\pi'}h'$, with $h\in H(\lambda)$
and $h'\in H(\mu)$. Then $E_\lambda(a)\omega_\pi E_\mu(b)=\phi_b(h)\phi_a(h') 
E_\lambda(a)\omega_{\pi'} E_\mu(b)$.

 We can thus replace 
$\pi$ by $\pi'$, and assume
that there are two adjacent cells, $i$ and $i+1$,
 with $\pi(i+1)=\pi(i)+1=j+1$, where $(i\;i+1)\in R(\lambda)$ 
and $(j\;j+1)\in R(\mu)$.
Since strings $i$ and $i+1$ finish at
adjacent points, any other string in $\omega_\pi$ either crosses above or below
both strings. Then
$\sigma_i\,\omega_\pi=\omega_\pi\,\sigma_j$, and so 
$(\sigma_i-a)\omega_\pi=\omega_\pi(\sigma_j-a)$.
Knowing that $(i\;i+1)\in R(\lambda)$  we can use Cor.~\ref{scalar} 
to write $E_\lambda(a)=h(\sigma_i -a)$, and similarly $E_\mu(b)=(\sigma_j-b)h'$.
Then \begin{eqnarray*}
E_\lambda(a)\omega_\pi E_\mu(b)&=& h(\sigma_i-a)\omega_\pi (\sigma_j-b)h'\\
&=&h\omega_\pi(\sigma_j-a)(\sigma_j-b)h'=0.
\end{eqnarray*}
 Interchange $a$ and $b$ to get the other equation.
\end{proof}
\begin{ascor}
\label{insep}
If $\lambda$ and $\mu$ are inseparable then 
\[
E_\lambda(a)H_n E_\mu(b)=0=E_\lambda(b)H_nE_\mu(a)\;.
\] 
\end{ascor}
\begin{proof} 
$H_n$ is spanned by $\{\omega_\pi;\pi\in S_n\}$, and no $\pi$ 
separates $\lambda $ from $\mu$.
\end{proof}
\begin{aslemma}
\label{justsep}
 Let $\lambda $ have $n$ cells. Then
$E_\lambda(a)\omega_\pi E_{\lambda^\vee}(b)$ is a scalar multiple of 
$E_\lambda(a)\omega_{\pi_\lambda} E_{\lambda^\vee}(b)$ for each $\pi\in S_n$.
\end{aslemma}
\begin{proof} 
By Lemma~\ref{=zero}, we need only consider the case when $\pi$
separates $\lambda$ from $\lambda^\vee$.

In this case $\omega_\pi= \omega_\rho\omega_{\pi_\lambda}\omega_\sigma$ with
$\rho\in R(\lambda), \sigma\in R(\lambda^\vee)$.  Then
\[
E_\lambda(a)\omega_{\pi} E_{\lambda^\vee}(b)=
E_\lambda(a)\omega_\rho\omega_{\pi_\lambda}\omega_\sigma E_{\lambda^\vee}(b)\\
=\phi_a(\omega_\rho)\phi_b(\omega_\sigma)
E_\lambda(a)\omega_{\pi_\lambda} E_{\lambda^\vee}(b).
\]
\end{proof}

\begin{astheorem}
\label{main}
Let $\lambda$ and $\mu$ be Young diagrams with $n$ cells.
Then 
\begin{eqnarray*}
e_\lambda e_\mu &=&0\quad\mbox{ for $\lambda\neq\mu$},\\
e_\lambda^2&=&\alpha_\lambda e_\lambda \quad \mbox{ for some scalar } \alpha_\lambda.
\end{eqnarray*}

Thus distinct
Young diagrams determine orthogonal elements,
 while each $e_\lambda$ is a quasi-idempotent element of $H_n$.
\end{astheorem}
\begin{proof}
\label{promai}
By definition, 
\[
e_\lambda e_\mu=E_\lambda(a)\omega_{\pi_\lambda}E_{\lambda^\vee}(b)
				{\omega}_{\pi_\lambda}^{-1}
E_\mu(a)\omega_{\pi_\mu}E_{\mu^\vee}(b){\omega}_{\pi_\mu}^{-1}\;.
\]
If $\lambda\neq\mu$ then either $\lambda^\vee$ and $\mu$ are inseparable, and so
$E_{\lambda^\vee}(b)\left(\omega_{\pi_\lambda}^{-1}\right)E_\mu(a)=0$, or $\lambda$ and
$\mu^\vee$ are inseparable and then $E_\lambda(a)\left(\omega_{\pi_\lambda}
E_{\lambda^\vee}(b){\omega}_{\pi_\lambda}^{-1}
E_\mu(a)\omega_{\pi_\mu}\right)E_{\mu^\vee}(b)=0$.
Thus $e_\lambda e_\mu=0$ when $\lambda\neq\mu$.

When $\lambda=\mu$ we can write 
\[
E_\lambda(a)\left(\omega_{\pi_\lambda}E_{\lambda^\vee}(b)
{\omega}_{\pi_\lambda}^{-1}
E_\lambda(a)\omega_{\pi_\lambda}\right)E_{\lambda^\vee}(b)=\alpha_\lambda 
E_\lambda(a)\omega_{\pi_\lambda}E_{\lambda^\vee}(b), \]
by Lemma~\ref{justsep}.
Then $e_\lambda^2=\alpha_\lambda 
E_\lambda(a)\omega_{\pi_\lambda}E_{\lambda^\vee}(b)\omega_{\pi_\lambda}^{-1}=
\alpha_\lambda e_\lambda$.
\end{proof}

The scalar $\alpha_\lambda$ \label{alpha} is calculated explicitly in 
\cite{qdim} as
\[
\alpha_\lambda=\prod_{(i,j)\in\lambda} s^{j-i}[\lambda_i+\lambda^\vee_j-i-j+1]\;,
\]
where $[k]=\displaystyle{s^k-s^{-k}\over s-s^{-1}}$ denotes the `quantum integer $k$'.

Note that $\lambda_i+\lambda^\vee_j-i-j+1$ is the {\it hook length\/} of the 
cell in the $i$th row and $j$th column.

%% file: AMsec5
\section{Framing and the Murphy operators}

In this section we return to the $3$-dimensional views of the skein 
${\cal S}(B,P)$ for a ball $B\cong B^3$ with a set $P$ of $2n$ distinguished
points on its boundary, introduced in  Sect.~\ref{idemp}.

We noted earlier that this skein is linearly isomorphic to $H_n$. In order to
translate the results about $H_n$ from the previous section into a
$3$-dimensional setting we shall use some explicit choices of isomorphism in a
couple of standard cases.

\begin{figure}
\[
\epsfxsize.9in\epsffile{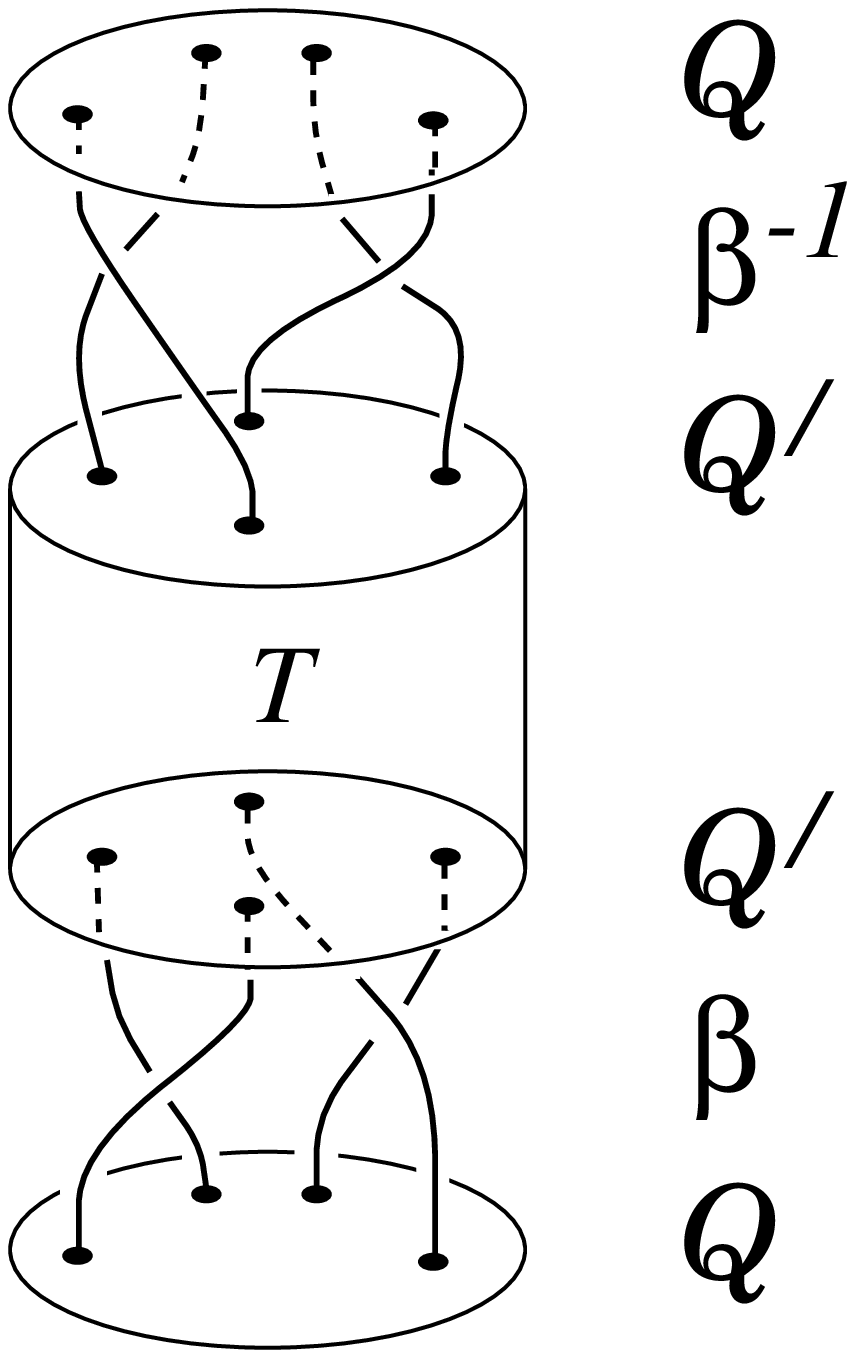}
\]
\caption{The wiring of $(D^2\times I,P^\prime)$ into $(D^2\times I,P)$}
\label{pprime}
\end{figure}

We start with the case when $B=D^2\times I$ and the set $P=P_I\cup P_O$ of
input and output points consists of $P_I=Q\times \{1\}$ and $P_O=Q\times\{0\}$
for some $Q\subset D^2$ with $|Q|=n$. Write  ${\cal S}(B,P)=H_Q$ in this case,
which will be identified with $H_n$ when the points of $Q$ lie in a straight
line across $D^2$. We can compare
$H_{Q'}$ and $H_{Q}$ by choosing  an $n$-string braid
$\beta$ in $D^2\times I$ with end points $Q\times\{0\}$ and $Q'\times\{1\}$, or
equivalently a homeomorphism from $(D^2,Q)$ to $(D^2,Q')$. Now construct a
wiring as shown in Fig.~\ref{pprime}, where the strings meet the level disks in $Q$ or
$Q'$ as indicated and are connected by $\beta$ or $\beta^{-1}$. This induces a
linear map ${\cal S}(\beta):H_{Q'}\to H_{Q}$, where the image
of a tangle $T$ in $H_{Q'}$ is the extended tangle in $H_Q$
shown in Fig.~\ref{pprime}.

 Note that $\beta$ is a braid in $D^2\times I$ in the sense that its
strings are monotone in the last coordinate. It should properly be regarded as
an element of the `braid groupoid' of $D^2$.

The skein $H_Q$ is clearly an algebra under the obvious stacking operation, and
${\cal S}(\beta)$ is an algebra homomorphism. Indeed it is an algebra 
isomorphism with inverse  ${\cal S}(\beta^{-1})$. Thus $H_Q\cong H_n$ for each
$Q$ with $|Q|=n$.

In the previous section we constructed an element $E_\lambda\in H_Q$ where the
points $Q$ lie in the cells of the Young diagram $\lambda$. The element
$e_\lambda\in H_n$ has the form $e_\lambda={\cal S}(\beta)(E_\lambda)$ where
the braid  $\beta$  lines up the cells of $\lambda$ according to the tableau
$T(\lambda)$. Then $E_\lambda$ is itself a quasi-idempotent in $H_Q$.
Now $E_\lambda$  is constructed as the composite $E_\lambda=E^I_\lambda
E^O_\lambda$  of  two  elements of
$H_Q$, with
$E^I_\lambda$ made from the row boxes with elements
$a_j$ and
$E^O_\lambda$ made from the column boxes with $b_j$.
 In the notation of the last
chapter we have
${\cal
S}(\beta)(E^I_\lambda)=E_\lambda(a)$ and ${\cal
S}(\beta)(E^O_\lambda)=\omega_{\pi_\lambda}
E_{\lambda^\vee}(b)\omega_{\pi_\lambda}^{-1}$. To understand the behaviour of
$e_\lambda $   further, by considering $E_\lambda$ and other elements of the
skeins $H_Q$,  we shall show the following.

\begin{aslemma}
\label{onedim} Let $Q\subset  D^2$ with $|Q|=n$. There is a unique element
$a_Q$ in the algebra $H_Q$ which
corresponds to $a_n\in H_n$ under {\em every} choice of isomorphism  ${\cal
S}(\beta)$.
\end{aslemma}

\begin{proof} The element $a_n\in H_n$ is central, by Lemma~\ref{linhom}.
Given $Q$ select any braid $\beta$ to connect $Q$ to the straight
line of $n$ points, and define $a_Q={\cal S}(\beta)(a_n)$. Let $\gamma$ be 
another such braid. Then $\gamma^{-1}\beta$ represents an element of $H_n$ which
commutes with  $a_n$.  Hence 
${\cal S}(\gamma^{-1}\beta)(a_n)=a_n$, and ${\cal S}(\gamma)(a_n)={\cal
S}(\beta)(a_n)$. 

A similar element $b_Q$ can be constructed from $b_n$.
\end{proof}

For every subdisk $D'\subset D^2$ with $R=Q\cap D'$ there is an induced
inclusion of algebras $H_R\subset H_Q$, where a tangle in $D'\times I$ is
extended by the trivial strings on $Q - R$. The discussion following
Lemma~\ref{linhom} shows that $a_Qa_R$ is a non-zero multiple of $a_Q$, and
similarly $b_Qb_R$ is a non-zero multiple of $b_Q$.

  We now look at some
consequences of the work in Sect.~4 within the context of the general skein
${\cal S}(B,P)$.  Take
$B=B^3$ and
$P=P_O\cup P_I\subset S^2$. Define a {\it geometric partition\/} $\omega$ of $P$
to be a family of disjoint discs $\{D_\alpha\}$ in $S^2$ containing the points
of $P$, such that no disc $D_\alpha$ contains both output and input points.

The geometric partition $\omega$ determines two partitions $\lambda(\omega)$
and $\mu(\omega)$ of $n$, where $\lambda_1\ge\lambda_2\ge\cdots\ge\lambda_k$
are the numbers of points of $P_I$ in the individual disks $D_\alpha$ of the
partitioning family, and $\mu$ is determined similarly by the output points
$P_O$.

 Given a geometric partition $\omega$ construct a wiring in $S^2\times I$ as
follows. For each disk $D_\alpha$ containing a subset $P_\alpha\subset P$
insert the skein element $a_{P_\alpha}$ or $b_{P_\alpha}$ into $D_\alpha\times
I\subset S^2\times I$, choosing $a_{P_\alpha}$ if $P_\alpha\subset P_O$ and
$b_{P_\alpha}$ if $P_\alpha\subset P_I$. The union of these gives a skein element
in $S^2\times I$, which induces a linear map $${\cal S}(\omega):{\cal S}(B,P)\to
{\cal S}(B,P)$$ by attaching $S^2\times I$ as a `shell' around $B^3$.

In many cases the nature of the map ${\cal S}(\omega)$ depends only on the
partitions $\lambda(\omega)$ and $\mu(\omega)$, as described in the following
lemma, which is an immediate consequence of Lemmas~\ref{=zero} and
 \ref{justsep}.

\begin{aslemma}
\label{lute}

Let ${\cal S}(\omega):{\cal S}(B,P)\to
{\cal S}(B,P)$ be the linear map induced from a geometric partition $\omega$ of
$P$. 

(a)\quad  If $\lambda(\omega)$ and $\mu(\omega)$ are inseparable  then ${\cal
S}(\omega)=0$.

(b)\quad If $\lambda(\omega)$ is just separable from $\mu(\omega)$ (when
$\mu(\omega)=\lambda(\omega)^\vee$) then ${\cal S}(\omega)$ has rank $1$. Its
image is spanned by ${\cal S}(\omega)(T)$, where $T$ is any tangle whose arcs
separate
$\lambda(\omega)$ from
$\mu(\omega)$.
\end{aslemma}
\begin{proof} Choose a homeomorphism from $B$ to $D^2\times I$ which carries
the disks $D_\alpha$ containing $P_O$ to $D^2\times\{0\}$ and those containing
$P_I$ to
$D^2\times \{1\}$, arranged in partition order along a straight line with the
input or output points lined up inside each disk. This induces a linear
isomorphism, $\phi$ say, from ${\cal S}(B,P)$ to $H_n$. The homeomorphism also
extends to a map of $S^2\times I$ which carries the wiring determined by
$\omega$ to the element $E_{\lambda(\omega)}(a)$ at the top of $D^2\times I$
along with
$E_{\mu(\omega)}(b)$ at the bottom. Then $\phi\circ{\cal
S}(\omega)\circ\phi^{-1}:H_n\to H_n$ is the linear map given by
$T\mapsto E_{\lambda(\omega)}(a)\; T\; E_{\mu(\omega)}(b)$. 

In case (a) this map is zero by Lemma~\ref{=zero}, and so    ${\cal
S}(\omega)=0$.

In case (b) its image has dimension 1, spanned by the image of any $T$ which
separates $\lambda(\omega)$ from
$\mu(\omega)$, by
Lemmas~\ref{pimult} and \ref{justsep}. The isomorphism $\phi$ gives the corresponding result for
${\cal S}(\omega)$ using tangles in ${\cal S}(B,P)$ which correspond to
separating tangles $T$ under the homeomorphism.
\end{proof}

Even where $\omega$ is separable we can deduce easily that ${\cal
S}(\omega)(T)=0$ for certain tangles $T$.

\begin{aslemma}
\label{subpart}
 Let $\omega'$ be a subpartition of $\omega$, in the sense that
every disk of $\omega'$ is contained in a disk of $\omega$. Then ${\cal
S}(\omega)\left({\cal
S}(\omega')(T)\right)$ is a non-zero multiple of ${\cal
S}(\omega)(T)$ for every $T$.
\end{aslemma}
\begin{proof}
This follows at once from the construction of ${\cal S}(\omega)$ and the
properties of the elements $a_Q$ and $b_Q$, when the shell for $\omega'$ is
viewed as lying inside the shell for $\omega$.
\end{proof}
Hence if ${\cal
S}(\omega')(T)=0$ then ${\cal
S}(\omega)(T)=0$ also.

\begin{aslemma}
\label{subball}
 Suppose that some of the strings of a tangle $T$ can be
enclosed in a ball $B'\subset B$ which does not meet the remaining strings of
$T$. A given geometric partition $\omega$ of $P=T\cap\partial B$ determines a
partition $\omega'$ of $P'=P\cap B'$, using  disks $D_\alpha\cap B'$. If
$\lambda(\omega')$ and $\mu(\omega')$ are inseparable then   ${\cal
S}(\omega)(T)=0$.
\end{aslemma}
\begin{proof}
 Extend $\omega'$ to a partition of $P$ by choosing the trivial partition on
the  points of $P-P'$. Then  ${\cal
S}(\omega')(T)=0$, by Lemma~\ref{lute} applied to $B'$. Now $\omega'$ is a
subpartition of $\omega$, so  ${\cal
S}(\omega)(T)=0$, by Lemma~\ref{subpart}.
\end{proof}

The simplest case occurs when $B'$ contains just two strings whose inputs and
outputs each lie in a disk in $B'$ contained in one of the disks of $\omega$.
Then $\lambda(\omega')=\mu(\omega')=(2)$ are inseparable partitions.

We return now to the study of $e_\lambda\in H_n$ by means of
$E_\lambda=E^I_\lambda E^O_\lambda$ in the skein $H_Q$ based on the Young
diagram $\lambda$.  Write $\omega$ for the geometric
partition grouping the input points by the rows of $\lambda$ and the output
points by its columns.  Then $\lambda(\omega)=\lambda$ and
$\mu(\omega)=\lambda^\vee$ are just separable, and ${\cal
S}(\omega)(T)=E^I_\lambda T E^O_\lambda$ for any
$T\in H_Q$. By Lemma~\ref{lute}(b) we can
then write $E^I_\lambda T
E^O_\lambda=\varphi_T E^I_\lambda 
E^O_\lambda=\varphi_T E_\lambda$ for some scalar $\varphi_T$.
In particular, if $T$ is  central then $T E_\lambda=\varphi_T E_\lambda$, and
consequently ${\cal S}(\beta)(T) e_\lambda=\varphi_T e_\lambda$ in $H_n$, for
the appropriate $\beta$.

Thus any central element $c\in H_n$ satisfies $ce_\lambda=c_\lambda e_\lambda$
for some scalar $c_\lambda$. The eigenvalue $c_\lambda$ may be found as
$\varphi_T$ when we choose $T={\cal S}(\beta^{-1})(c)$ in the appropriate 
skein $H_Q$. We will now calculate these
eigenvalues for two important central elements of $H_n$.

\begin{figure}
\[
F_n\qquad=\qquad\raisebox{-1cm}{\epsfxsize.5in\epsffile{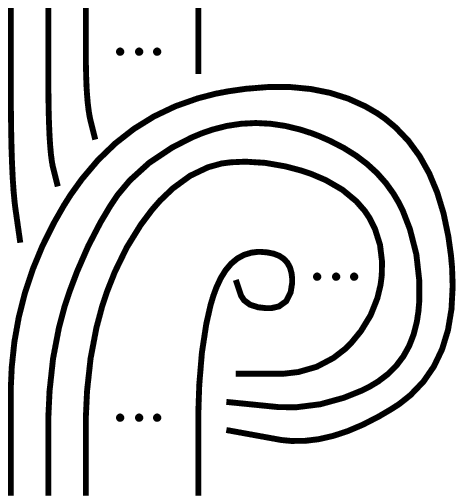}}
\]
\caption{The full curl $F_n$ on $n$ strings}
\label{fullcurl}
\end{figure}

The first of these is the full curl $F_n\in H_n$, which is closely involved
with changes of framing. In the skein $H_Q$ we can describe the full curl $F_Q$
 as follows. Take the trivial framed braid $Q\times I\subset D^2\times I$ and
 rotate $D^2\times\{1\}$ through $2\pi$ anticlockwise while keeping
$D^2\times\{0\}$ fixed. The resulting framed braid represents $F_Q$.
Note that  the framing of each string in $F_Q$ is $1$. It is clear that ${\cal
S}(\beta)F_{Q'}=F_Q$ for any $\beta$. The full curl $F_n$ can be drawn, with
blackboard framing, as in 
Fig.~\ref{fullcurl}.  We can describe $F_Q$ recursively by selecting a subdisk
containing $n-1$ of the points, $R$ say. Then $F_Q$  is the
composite of $F_R$ with the braid which takes the remaining string once around
all the others in the subdisk.

\begin{figure}
\[
\epsfxsize.6in\epsffile{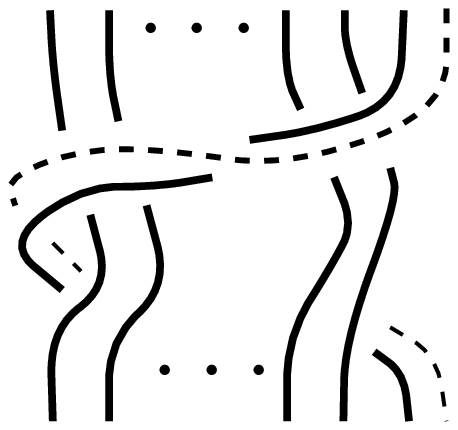}
\]
\caption{The braid $T_i$}
\label{ti}
\end{figure}

Using tensor product to denote the juxtaposition of oriented tangles in $H_n$,
the recursive definition gives 
\[
F_n=T_n(F_{n-1}\otimes 1)=\prod_{i=1}^n T_i\otimes 1^{n-i}\;,
\]
where $T_i$ is the $i$-string tangle in  Fig.~\ref{ti} (all strings 
have blackboard framing except the $i$th which has its framing indicated by 
the dotted line in the figure).

\begin{astheorem}
\label{flamb}

Let $\lambda$ be a Young diagram
with $n$ cells. Then $F_ne_\lambda=f_\lambda e_\lambda$, where 
$f_\lambda=x^{\vert\lambda\vert^2}v^{-\vert\lambda\vert}s^{n_\lambda}$
and  $n_\lambda=\displaystyle\sum_{(i,j)\in\lambda}2(j-i)$.
\end{astheorem}
\begin{proof}
By induction on the number of strings.
When $n=1$ we have $\lambda=\Box$. Then $e_\lambda$ is the single string, and
 $f_\lambda=xv^{-1}$ by the skein relations.  

Assume now that the result holds for all Young diagrams with fewer than $n$
cells.   Let $\lambda$  be a Young diagram with $n$
cells and let $\mu$ denote the Young diagram obtained from
$\lambda$ by deleting a fixed choice of extreme cell, with coordinates $(p,r)$.

We shall work in the skein $H_Q$ based on $\lambda$, and show that
$F_QE_\lambda=f_\lambda E_\lambda$ by calculating ${\cal S}(\omega) (F_Q)$.
Take $R\subset Q$ to be defined by the subdiagram $\mu$. We can include the
algebra $H_R$ in $H_Q$ by adjoining an extra trivial string to the skein. Then
$F_Q =F_R T$, where $T$ is a version of Fig.~\ref{ti} in which the extreme
string winds once around the strings from $R$. It is geometrically clear that
$T$ commutes with the subalgebra $H_R$.

Now since $\mu$ is a subdiagram of $\lambda$ we have $E^I_\lambda E^I_\mu=
kE^I_\lambda$ and  $E^O_\mu E^O_\lambda=lE^O_\lambda$ 
for some non-zero $k,l$. Then $F_Q E_\lambda=E_\lambda^I F_Q E_\lambda^O$
since $F_Q$ is central and
\begin{eqnarray*}
{\cal S}(\omega) (F_Q)=
E^I_\lambda F_Q
E^O_\lambda&=&(lk)^{-1}E^I_\lambda E^I_\mu F_R T E^O_\mu E^O_\lambda \\
&=&(lk)^{-1}E^I_\lambda (E^I_\mu F_R  E^O_\mu) T E^O_\lambda \\
&=&f_\mu (lk)^{-1}E^I_\lambda (E^I_\mu E^O_\mu) T E^O_\lambda \\
&=&f_\mu 
E^I_\lambda T
E^O_\lambda =f_\mu{\cal
S}(\omega) (T).\\
\end{eqnarray*}
  It is  enough to prove that 
\begin{equation}
\label{gordon}
E^I_\lambda T
E^O_\lambda=x^{2\vert\mu\vert+1}v^{-1}
						s^{2(r-p)}\, E_\lambda\;,
\end{equation}
since  then $F_QE_\lambda=f_\mu x^{2\vert\mu\vert+1}v^{-1}
						s^{2(r-p)}\, E_\lambda$. 
Rewriting $f_\mu$ by the induction hypothesis we obtain
\[
f_\lambda =x^{\vert\mu\vert^2+2\vert\mu\vert+1}
					v^{-1-\vert\mu\vert}s^{2(r-p)+n_\mu}
	 =x^{\vert\lambda\vert^2}v^{-\vert\lambda\vert}s^{n_\lambda}\;.
\]
\begin{figure}
\[
\raisebox{-1.5cm}{\epsfxsize2in\epsffile{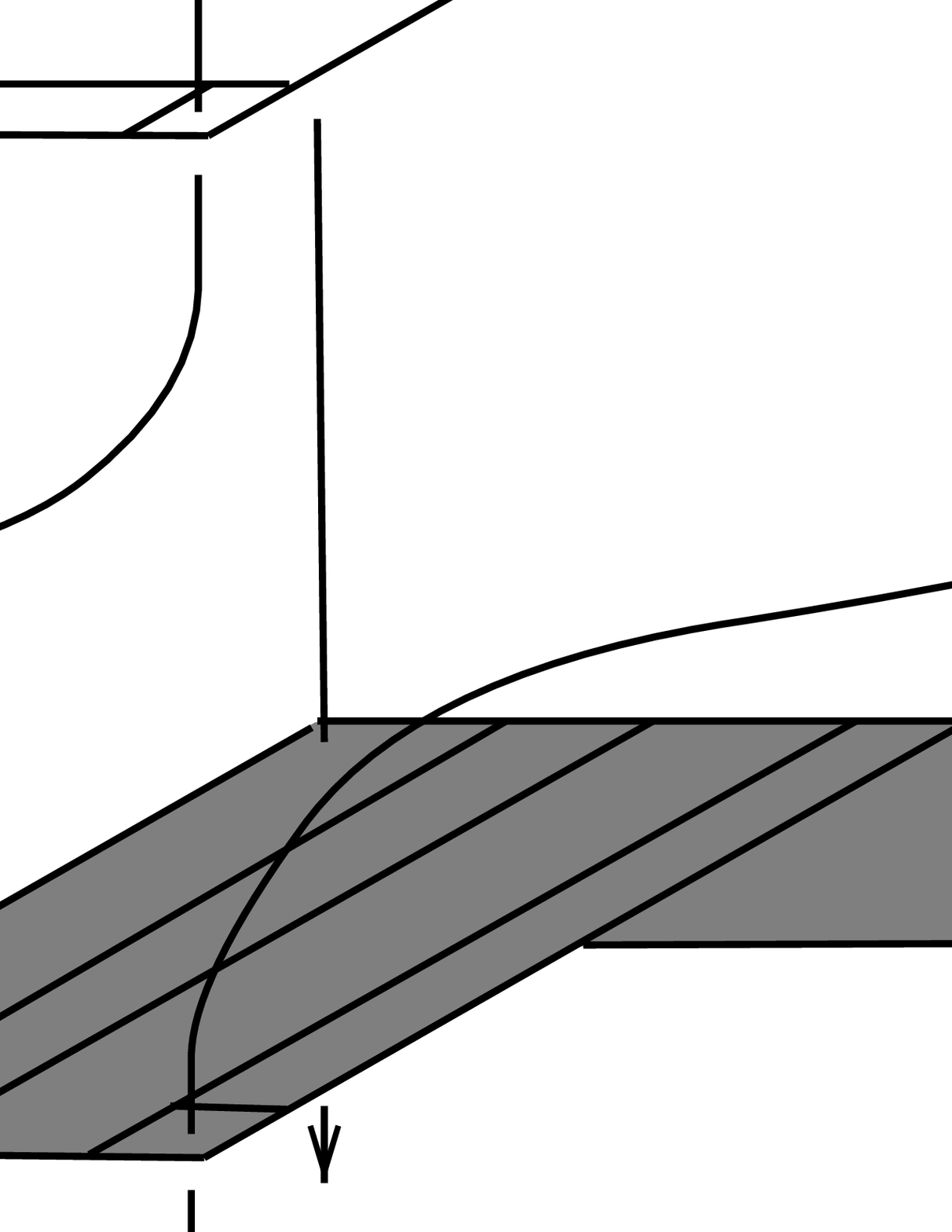}}
\]
\caption[]{The element $E^I_\lambda TE^O_\lambda={\cal S}(\omega)(T)$}
\label{richard}
\end{figure}
The  picture of $E^I_\lambda T
E^O_\lambda={\cal
S}(\omega) (T)$
is given in Fig.~\ref{richard}.

We  now represent some simple braids $B$ in the skein $H_Q$ by
schematic pictures looking down on $B$ from above.  The points of $Q$ are
points on the grid.  Any string not explicitly shown is assumed to pass
straight down, while a path with an arrow indicates the projection of a  string
descending from the top to the bottom in $D^2\times I$. When all strings pass
straight down we have the  identity braid $\mbox{Id}$, with
${\cal S}(\omega)(\mbox{Id})=E_\lambda$. 
  
The schematic picture for $T$  is shown  in Fig.~\ref{anna}. 
There is one string in this braid starting and
finishing from the chosen extreme position, marked  
$\times$, which encircles the vertical strings based on $\mu$. We now alter $T$
systematically through a sequence of braids which retain the vertical strings
while taking the encircling string round successively fewer verticals. When
we pass from one such braid $T'$ to another $T''$ by crossing through one
vertical string, noted on the accompanying diagrams by  
$\bullet$, the skein relation  gives $T'=x^2T''+x(s-s^{-1})L$, where the braid $L$
arises by smoothing at the crossing of the vertical and encircling strings. Thus
the diagram for $L$ has two non-vertical arcs, one from $\times$ to $\bullet$  and
one from $\bullet$  to $\times$. The braids $L$ which arise in this way will be
shown, using Lemma~\ref{subball}, to satisfy ${\cal S}(\omega)(L)=0$, giving
${\cal S}(\omega)(T')=x^2{\cal S}(\omega)(T'')$, when the string $\bullet$ does
not lie on the edge of the extreme rectangle determined by $\times$, in the
sense of Sect.~\ref{bib}.

Starting then with $T$, pull the encircling string successively through the
vertical strings outside the extreme rectangle for $\times$, working leftward
through the columns and upwards through the rows until it encircles the extreme
rectangle.
\begin{figure}
\[
T\quad=\quad\raisebox{-1cm}{\epsfxsize1in\epsffile{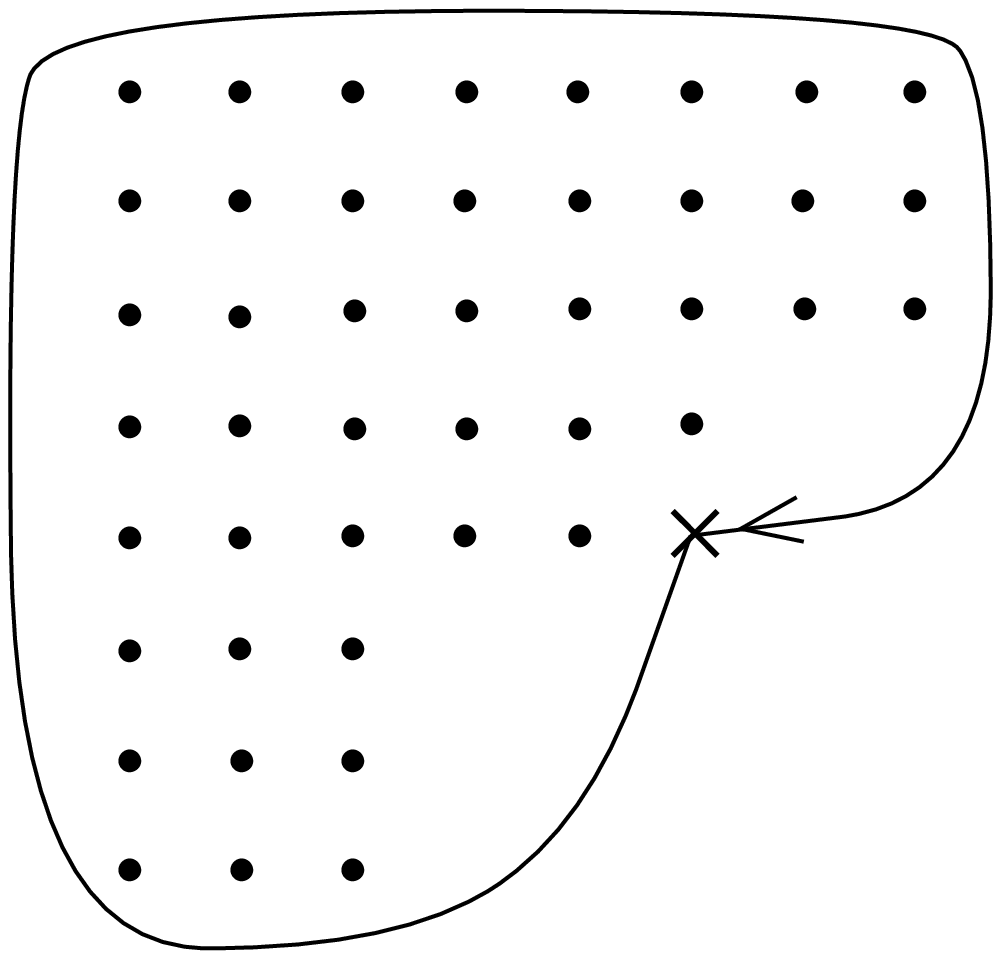}}\ ,\quad
T''\quad=\quad
\raisebox{-1cm}{\epsfxsize1in\epsffile{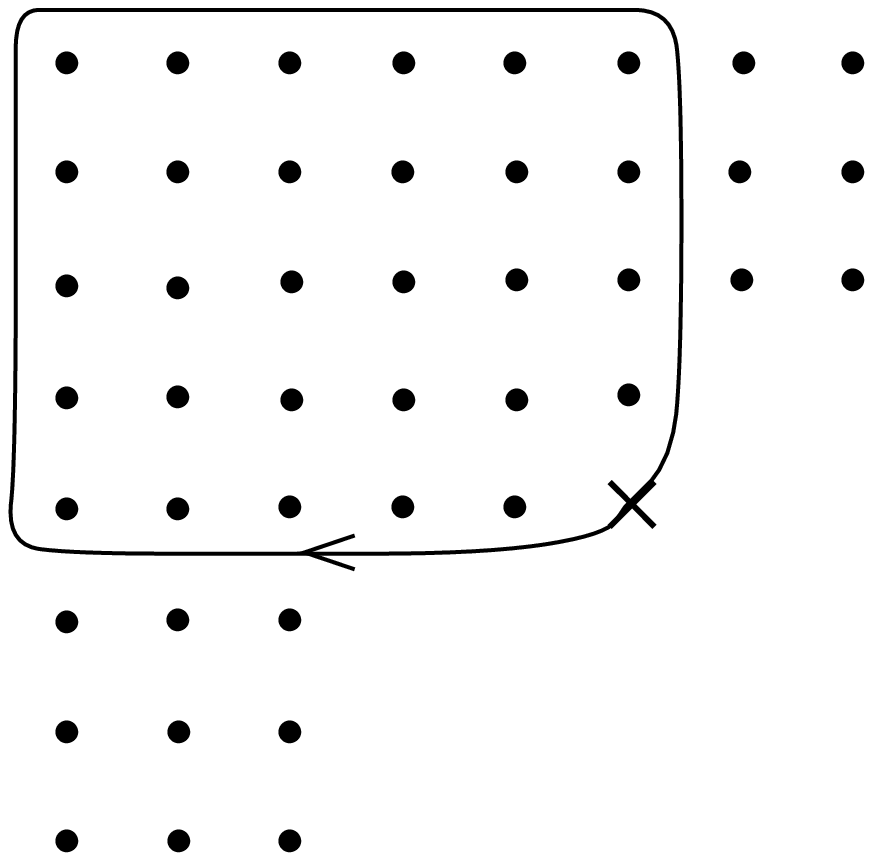}}
\]
\caption[]{Reduction to the extreme rectangle}
\label{anna}
\end{figure}
The diagram of a typical smoothed braid
$L$ is shown in  Fig.~\ref{dotpic}.
\begin{figure}
\[
\epsfxsize1in\epsffile{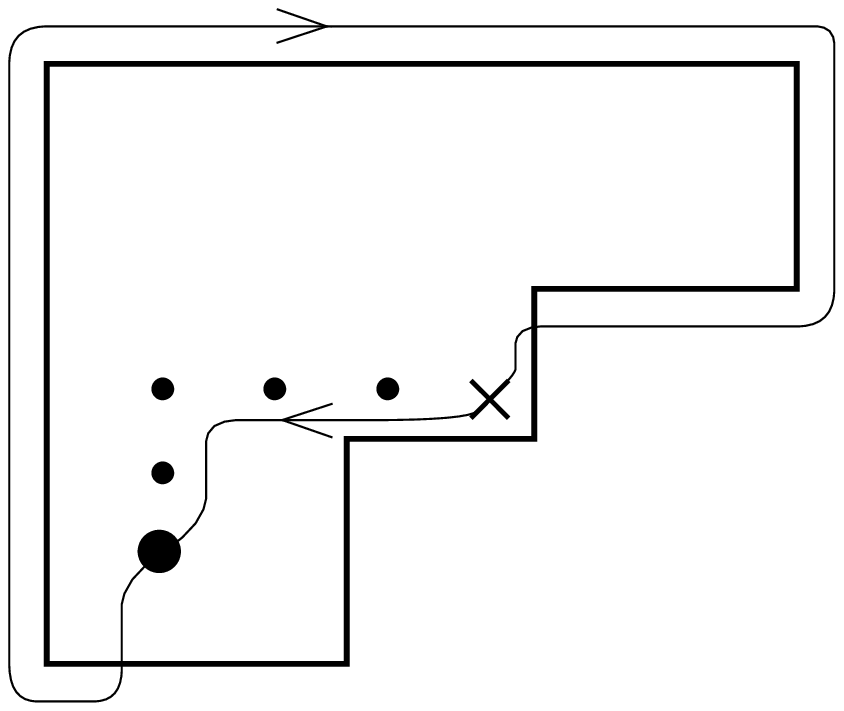}\phantom{this is a spacing device}\epsfxsize1in\epsffile{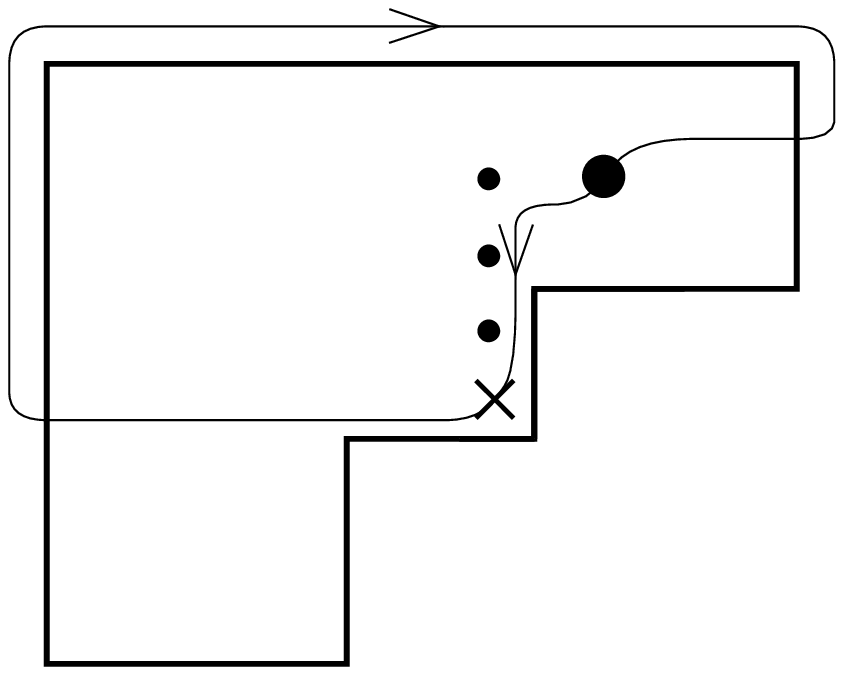}
\]
\caption{Smoothed braids arising from crossings outside the extreme rectangle}
\label{dotpic}
\end{figure}

Each $L$ has a
subbraid made up of an L-shaped array of vertical strings taken from the row
and column containing $\times$ and $\bullet$, along with a single string joining
$\times$ and $\bullet$. This subbraid can be enclosed in a ball $B'$ for which
the induced geometric partition $\omega'$ determines an inseparable pair of
partitions $\lambda(\omega')$ and $\mu(\omega')$. Thus ${\cal S}(\omega)(L)=0$
in each case, by Lemma~\ref{subball}. 
A typical example is shown in Fig.~\ref{lshape}.

\begin{figure}
\[
\epsfxsize.75in\epsffile{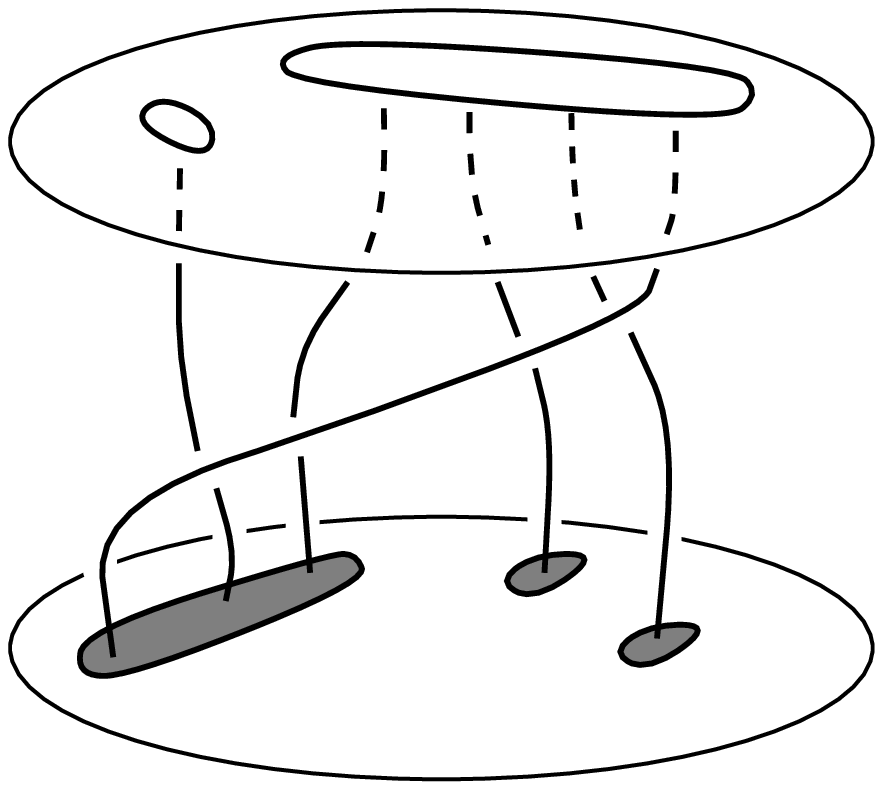}
\]
\caption{}
\label{lshape}
\end{figure}

The encircling string in $T$ can then be pulled through $|\lambda|-rp$ strings
to give the braid $T''$ in which it encircles the extreme rectangle, as in
Fig.~\ref{anna}, with
${\cal S}(\omega)(T)=x^{2(|\lambda|-rp)}{\cal S}(\omega)(T'')$.

Now write $T''=C_1 S C_2$, where $C_1, S $ and $C_2$ are shown in
Fig.~\ref{rectedge}.
\begin{figure}
\[
C_1\ =\ \raisebox{-1.3cm}{\epsfxsize1in\epsffile{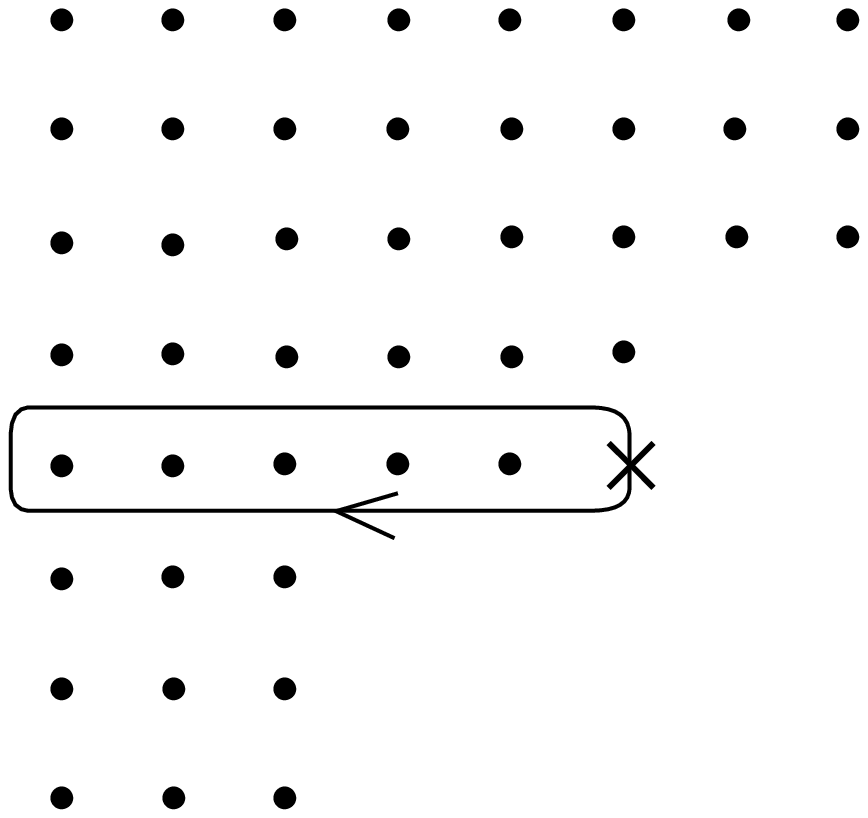}}\ ,\quad 
S\ =\ \raisebox{-1.3cm}{\epsfxsize1in\epsffile{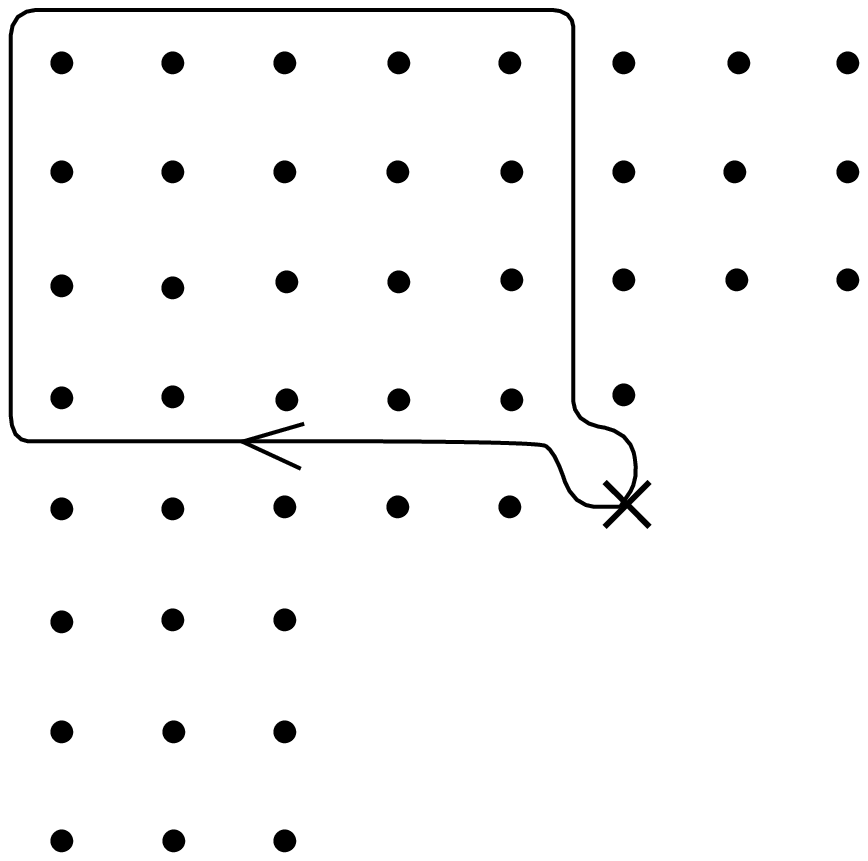}}\ ,\quad
C_2\ =\ \raisebox{-1.3cm}{\epsfxsize1in\epsffile{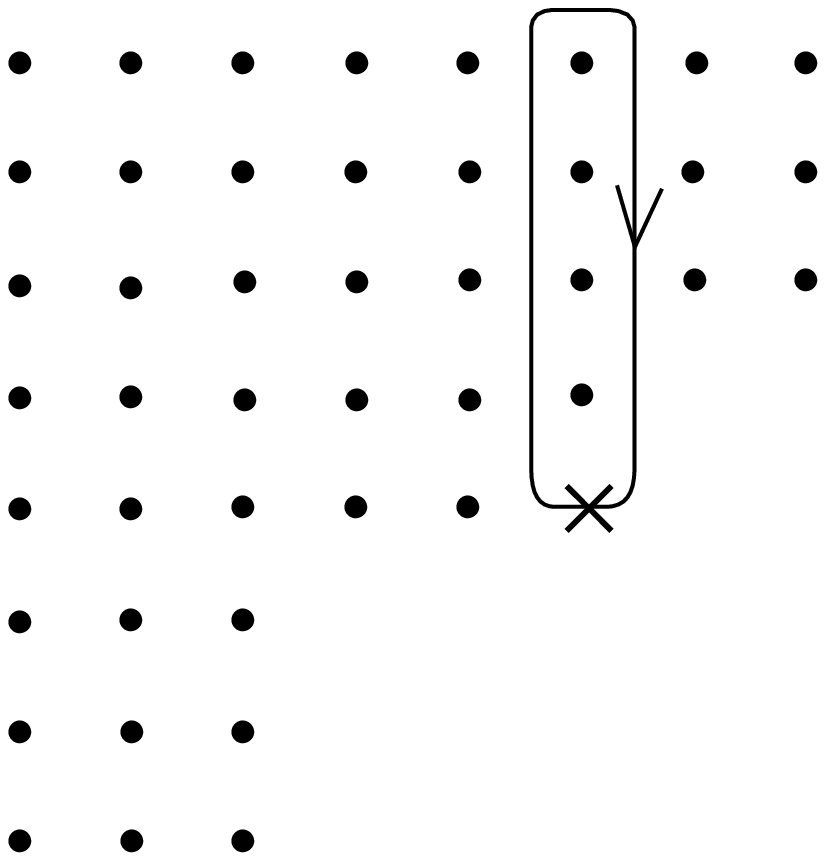}} 
\]
\caption{}
\label{rectedge}
\end{figure}
A similar systematic move of the encircling string in $S$ through the
vertical strings inside the extreme rectangle, working through the columns from
the left and from the top down, will replace $S$ by a braid with all strings
vertical, at the expense of a factor $x^{2(r-1)(p-1)}$. A typical smoothed braid
$L$ arising from the skein relation at an intermediate crossing change is shown
in Fig.~\ref{inex}, again with an L-shaped subdiagram  which ensures that $ {\cal
S}(\omega)(L)=0$.
\begin{figure}
\[
\epsfxsize1in\epsffile{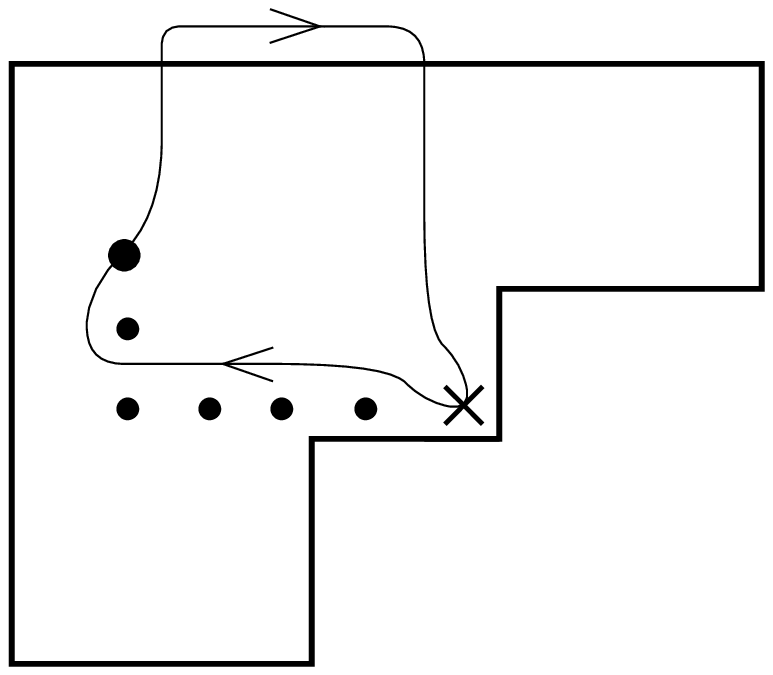}
\]
\caption{A smoothed braid arising from a crossing inside the extreme rectangle}
\label{inex}
\end{figure}

The framing of the eventual vertical string through
$\times$ is still
$+1$; when this is returned to
$0$ we get ${\cal S}(\omega)(S)=x^{2(r-1)(p-1)}xv^{-1}E_\lambda$. The
calculations for Eq.~(\ref{gordon}) are completed by noting that $E^I_\lambda
C_1=(xs)^{2(r-1)}E^I_\lambda$ and $C_2
E^O_\lambda=(-xs^{-1})^{2(p-1)}E^O_\lambda$, by Lemma~\ref{linhom}. 

Then ${\cal S}(\omega)(T'')=x^{2(r+p-2)}s^{2(r-p)}{\cal S}(\omega)(S)$, 
giving 
\[
{\cal S}(\omega)(T)=x^{2|\lambda|-1}s^{2(r-p)}v^{-1}E_\lambda
\]
as claimed.
\end{proof}

Note that
$\displaystyle\sum_{(i,j)\in\lambda}2j
	=\displaystyle\sum_{i=1}^k\displaystyle\sum_{j=1}^{\lambda_i} 2j
	=\displaystyle\sum_{i=1}^k \lambda_i(\lambda_i+1)
	=\displaystyle\sum_{i=1}^k\lambda_i^2+\vert\lambda\vert$.
Similarly,
$\displaystyle\sum_{(i,j)\in\lambda}2i=
\displaystyle\sum_{j=1}^m(\lambda_j^\vee)^2+\vert\lambda\vert$.
Thus, by Theorem~\ref{flamb}, 
$n_\lambda=\displaystyle\sum_i\lambda_i^2-\displaystyle\sum_j(\lambda_j^\vee)^2$.

As we discuss in Sect.~\ref{qinv}, Theorem~\ref{flamb} 
is a skein theoretic calculation of the framing factors for the irreducible 
representations of $SU(N)_q$.  
There are various other well-known formulae, in terms of the weights of the
representation or the length of either the columns or rows of the Young 
diagram.
In \cite{morton} Morton gives the $3$-variable formula for $f_\lambda$,
with $n_\lambda$ defined recursively in terms of smaller Young diagrams.
Alternatively, the fact that a Young diagram with exactly $N$ cells in the 
first column indexes the same representation as the Young diagram obtained 
by removing this column gives an inductive proof of a formula in 
terms of the column length.
For a given $N$, the framing factor becomes a function of $s$ alone.
Reshetikhin \cite{} proved that $f_\lambda(N)=s^{\Delta_\lambda}$ 
where $\Delta_\lambda$ is the value of the Casimir operator on the 
irreducible representation associated to $\lambda$. 
An expression for $\Delta_\lambda$ in terms of the weights of the 
representation, was described to the authors by Kohno, from
which a form of $n_\lambda$ in terms of row length alone can be derived,
namely, $n_\lambda=\sum\lambda_i^2+\vert\lambda\vert-\sum2i\lambda_i$.

However, this lacks the simplistic charm of stating that $n_\lambda$ is 
twice the sum of the contents.

The other central element of $H_n$ which we consider is
the sum $M$ of the Murphy operators.
These are discussed in \cite{dip} and are most easily described in the version
$H_n(1,z)$ of the  Hecke
algebras, converting by the appropriate isomorphism as needed.
In $H_n(1,z)$ the $j$th {\it Murphy operator\/} 
$M(j)=\displaystyle{\sum_{i=1}^{j-1}\omega_{(i\,j)}}$ is 
the sum of all the positive {\it transposition\/} braids involving the 
string $j$ and strings of
smaller label. Then write $M=\sum_{j=1}^n M(j)=\sum \omega_\tau$, where
the sum is over all transpositions  $\tau\in S_n$, to give an element which
can readily be shown by skein theory to be central in $H_n$. We know from our
work above that eigenvalues $m_\lambda$ can be found such that 
$Me_\lambda=m_\lambda e_\lambda$. 
As above we shall calculate $m_\lambda$ from
the equation  
$E^I_\lambda {\cal S}(\beta)(M)E^O_\lambda=m_\lambda E_\lambda$ 
in $H_Q$, with $Q$ defined by the diagram
$\lambda$. We again use schematic diagrams to consider the 
constituent braids  
${\cal S}(\beta)(\omega_{(i\,j)})=\beta^{-1}\omega_{(i\,j)}\beta$ inside
$H_Q$. This allows an easy calculation of the scalar $m^{(j)}_\lambda$ 
which satisfies 
$E^I_\lambda {\cal S}(\beta)(M(j))E^O_\lambda=m^{(j)}_\lambda E_\lambda$ and 
hence gives $m_\lambda=\sum_j m^{(j)}_\lambda$. It provides an alternative to the original algebraic proof of 
Theorem~\ref{dip} by Dipper and James \cite{dip}.

One of the advantages of the central element $M$ over, say, the
full curl is demonstrated in \cite{chak}. 
There are many examples of distinct Young diagrams with 
the same framing factor, which is the eigenvalue of the full curl.
  The smallest example occurs with six cells, where
the partitions $(4,1,1)$ and $(3,3)$ both have framing 
factor $x^{36}v^{-6}s^6$.  Any pair of Young diagrams with the 
same number of cells, which are
self conjugate will have the same framing factor, for example $(4,2,1,1)$ 
and $(3,3,2)$.  In contrast, the scalar $m_\lambda$
uniquely determines $\lambda$, given the number of cells.  

Recall that the Young tableau $T(\lambda)$ was defined by labelling the cells
of $\lambda$ from $1$ to $\vert\lambda\vert$ along the rows.
For the three-dimensional picture, we use this to indicate the ordering of the
strings within the Young diagram, and to determine the conversion braid $\beta$
used in the isomorphism
${\cal S}(\beta)$ with the Hecke algebra.

As we saw in the case of the framing coefficient, the power of $x$
really only keeps track of the number of strings.  We shall set $x=1$
in what follows. This gives one of the presentations
of the Hecke algebra, $H_n(1,z)$,
described in Sect.~\ref{iso}.  The result for $H_n$
follows directly from an application of the isomorphism with $H_n(1,z)$.

\begin{figure}
\[
\omega_{(i\,j)}\quad=\quad\raisebox{-3mm}{\epsfxsize1in\epsffile{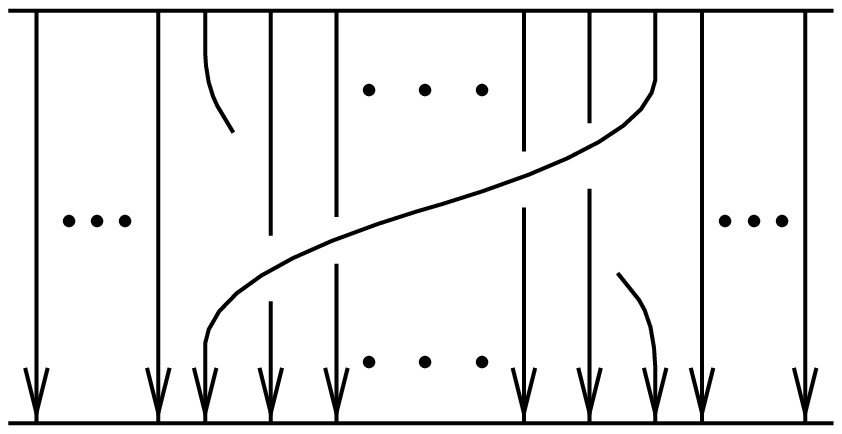}}
\]
\caption{}
\label{bij}
\end{figure}

We can write $M(j)$ as the sum of positive transposition braids 
\[
M(j)=\sum_{i=1}^j \omega_{(i\,j)}\;,
\]
where $\omega_{(i\,j)}$ is the positive permutation braid in Fig.~\ref{bij}.

\begin{aslemma}
\label{peter}
Suppose that the $j$th cell of $\lambda$ has coordinates $(k,l)$. Then
\[
E^I_\lambda {\cal S}(\beta)(M(j))E^O_\lambda
	=s^{l-k}[l-k]\,\,E_\lambda\;.
\]
\end{aslemma}
\begin{proof}
Work as above in the algebra $H_Q$ based on the diagram $\lambda$. Write $j@i$
for the braid drawn schematically in Fig.~\ref{roundi}.
\begin{figure}[h]
\[
j@i\ =\ \raisebox{-1cm}{\epsfxsize1in\epsffile{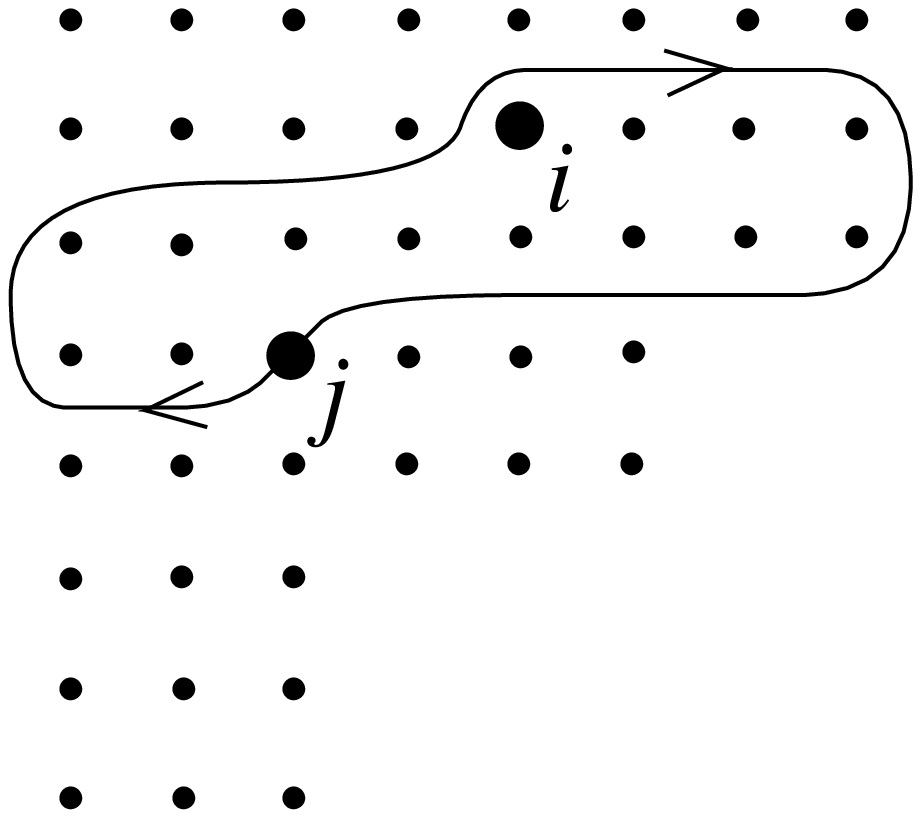}}\ ,\qquad
{\cal S}(\beta)(\omega_{(i\,j))})\ =\ 
\raisebox{-1cm}{\epsfxsize1in\epsffile{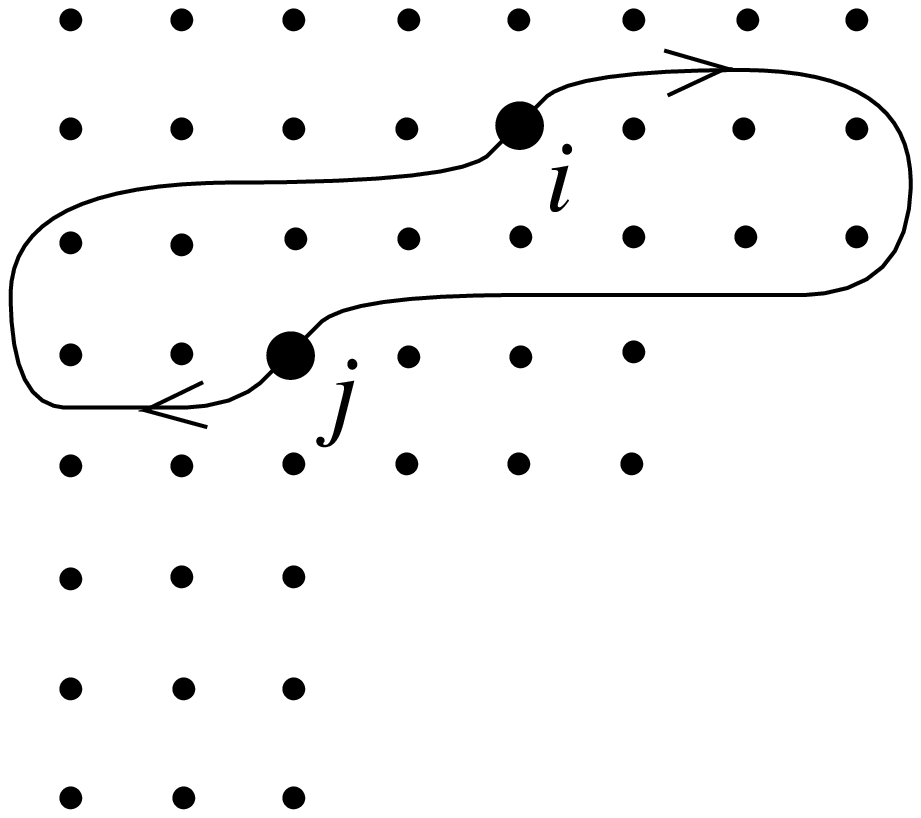}}
\]
\caption{}
\label{roundi}
\end{figure}
If we  pull the string $j$ through the vertical string $i$ we get the braid
$j@(i+1)$, while smoothing at the crossing of these strings gives the braid
${\cal S}(\beta)(\omega_{(i\,j)})$, shown also in Fig.~\ref{roundi}.
Then
\[
j@i=(s-s^{-1}){\cal S}(\beta)(\omega_{(i\,j)}) +j@(i+1),
\] by the  skein relation (with $x=1$).
Summing over all $i<j$ gives
\[
j@1-j@j=(s-s^{-1}){\cal S}(\beta)(M(j)).
\]
From the calculation  in Theorem~\ref{flamb} we have 
$${\cal S}(\omega)(j@1)=E_\lambda^I (j@i)E_\lambda^O=s^{2(l-k)}E_\lambda$$ while
${\cal S}(\omega)(j@j)=E_\lambda$, where $\omega$ is the geometric
partition determined, as before,  by the diagram $\lambda$.  Therefore
\[
E^I_\lambda {\cal S}(\beta)(M(j))E^O_\lambda={\cal
S}(\omega)(j@1-j@j)={(s^{2(l-k)}-1)\over (s-s^{-1})}\ \ E_\lambda\;.
\]
\end{proof}

We can thus prove the result of Dipper and James \cite{dip}
 directly using skein theory.
For completeness, we also include the proof that the eigenvalues of the 
sum of Murphy operators uniquely determines the Young diagrams, 
which is due to Katriel, Abdesselam and Chakrabarti \cite{chak}

\begin{astheorem}
\label{dip}

Let $\lambda$ be a Young diagram with $n$ cells.
Set $M=\sum_{j=1}^n M(j)$.
Then $M$ is central in $H_n$ and
$M e_\lambda=m_\lambda e_\lambda
$
where
\[
m_\lambda=\sum_{(k,l)\in\lambda}s^{(l-k)}[l-k] = \sum_{(k,l)\in\lambda}
{s^{2(l-k)}-1\over s-s^{-1}}\;.
\]
The Young diagram
$\lambda$ is determined by $m_\lambda$ and $n$.
\end{astheorem}
\begin{proof}

First we show that $M$ is central.  It is enough to show that
 $\sigma_i$ commutes with $M$, for $i=1,\ldots, n-1$.
We know that $\displaystyle{M=\sum_{k< l}\omega_{(k\,l)}}$ is the sum of
all the positive transposition braids, working as above in $H_n(1,z)$.
If neither $k$ nor $l$ is equal to $i$ or $i+1$ then 
$\sigma_i$ commutes with $\omega_{(k\,l)}$.
It also commutes with $\omega_{(i\;i+1)}=\sigma_i$. 

The remaining elements in $M$ can be written as
\[
\displaystyle{A_i=\sum_{k<i}\left(\omega_{(k\;i)}+\omega_{(k\;i+1)}\right)
+\sum_{i+1<l}\left(\omega_{(i\;l)}+\omega_{(i+1\;l)}\right)}\;.
\]
Now $\sigma_i\,\omega_{(k\;i)}\,\sigma_i=\omega_{(k\;i+1)}$, for $k<i$.
Similarly  $\sigma_i\,\omega_{(i+1\;l)}\,\sigma_i=\omega_{(i\;l)}$, for $i+1<l$.

Hence $A_i=X+\sigma_i\, X\,\sigma_i$, where $X=\displaystyle{\sum_{k<i}\omega_{(k\;i)}
+\sum_{i+1<l}\omega_{(i+1\;l)}}$.
Then
\begin{eqnarray*}
\sigma_i\, A_i&=&\sigma_i\,X+\sigma_i^2\,X\,\sigma_i\\
&=&\sigma_i\,X+z\,\sigma_i\,X\,\sigma_i+X\,\sigma_i\\
&=&A_i\,\sigma_i,
\end{eqnarray*}
since $\sigma_i^2=z\,\sigma_i+1$ in $H_n(1,z)$. 
This completes the proof that $M$ is central.

Then  ${\cal S}(\beta)(M)$ is central in $H_Q$ and so ${\cal
S}(\beta)(M)E_\lambda=E_\lambda^I {\cal S}(\beta)(M) E_\lambda^O =m_\lambda
E_\lambda$. The  value of $m_\lambda$ as claimed follows directly from
Lemma~\ref{peter} by writing $M=\sum_j M(j)$.

To work in the  version of the Hecke algebra used by Dipper and James we need
simply apply the  isomorphism from $H_n(1,z)$ to $H_n(q)$. A further small
change is needed  to give their exact result, as the Murphy operators in their
paper correspond to $s\,M(j)$ here, and so their value of $m_\lambda$ is $s$
times our value.

To recover $\lambda$ from $m_\lambda$ calculate the Laurent polynomial
\[ (s-s^{-1})m_\lambda+n=\sum_{(k,l)\in\lambda}s^{2(l-k)}=\sum c_js^{2j}
\mbox{ say.}\]
Then the integer $c_j$ gives the number of cells in $\lambda$ on the diagonal
at distance $j$ to the right of the leading diagonal.
\end{proof}

%% file: AMsec6
\section{Connections with the quantum group invariants}
\label{qinv}

Here we discuss briefly the relationship between the Homfly polynomial
and the $SU(N)_q$--invariants and the bearing the idempotent elements 
of Sect.~\ref{idemp} have on this relationship.
For more details, the reader is referred to \cite{mine} and \cite{knot96}.

Let $Q_\lambda$ denotes the closure 
of $(1/\alpha_\lambda)\, e_\lambda$ in the Homfly skein
of the annulus, or equally of $(1/\alpha_\lambda)\, E_\lambda$ in the skein of $D^2\times S^1$. 

\begin{astheorem}[\rm\protect\cite{tur}]
\label{homq}
The quantum invariant $J(L;V_\Box,\cdots, V_\Box)$ 
of the link $L$ with each component decorated by the fundamental representation
is given as a function of $s$ by the 
framed Homfly polynomial ${\cal X}(L)$, evaluated at
$x=s^{-1/N}$ and $v=s^{-N}$.
We will denote this evaluation of ${\cal X}$ by ${\cal X}_N$.
\end{astheorem}

Every element of ${\cal S}(R_n^n)$ determines
an endomorphism of the $n$-fold tensor product of
the fundamental representation, $V_\Box^{\otimes n}$.
Let $\Phi$ denote the map which takes an
$(n,n)$ tangle to this endomorphism.  
Jimbo \cite{jimbo} showed that $\Phi$ is a surjective ring homomorphism.

\begin{asprop}[\rm\protect\cite{mine}]
For every Young diagram $\lambda$ with $n$ cells, the endomorphism 
of $V_\Box^{\otimes n}$ determined by $\Phi(e_\lambda)$ is
a scalar multiple of the
projection map onto a single copy of the irreducible
$SU(N)_q$--module $V_\lambda$.
\end{asprop}
\begin{astheorem}[\rm\protect\cite{mine}]
\label{pork}
  
Let $C$ be a framed knot coloured by the irreducible representation
$V_\lambda$.
Let $S$ be the satellite knot $C*Q_\lambda$
with companion $C$
and pattern $Q_\lambda$.
Then
\[
J(C;V_\lambda)={\cal X}_N(S)\;.
\]
The result also holds for links where each component coloured
by $V_\lambda$ is replaced by its satellite, with pattern $Q_\lambda$.
\end{astheorem}

It is not difficult to show that the elements $Q_\lambda$ form a linear
basis for the skein of the annulus, by expressing them in terms
of Turaev's basis, as discussed in Sect.~\ref{bazza}.  There is no
nice skein theory proof that the expression for $Q_\lambda Q_\mu$
in terms of this basis is identical to that for the decomposition of
the tensor product of the irreducible representations $V_\lambda\otimes V_\mu$.
It can, however, be proved more circuitously that
\[
Q_\lambda Q_\mu=\sum_{\nu\in\lambda\mu} a_{\lambda\mu\nu} Q_\nu\;,
\]
where the $a_{\lambda\mu\nu}$ are the Littlewood-Richardson coefficients.

Theorem~\ref{pork} has been used 
to calculate the quantum dimensions of representations of $SU(N)_q$.
The details are discussed in \cite{qdim}, where the following formula 
is proved skein theoretically;
\[
{\cal X}(\widehat{e}_\lambda)= \prod_{(i,j)\in\lambda} s^{j-i}
			{(v^{-1}s^{j-i}-vs^{i-j})\over(s-s^{-1})}\;.
\]
Making the substitutions prescribed by Theorem~\ref{homq}, we obtain the
 formula for the quantum dimension ${\cal X}_N(Q_\lambda)$, which
was originally established by Reshetikhin \cite{r1},
\[
{\cal X}_N(Q_\lambda)=
\prod_{(i,j)\in\lambda} {[N+j-i]\over [\lambda_i+\lambda^\vee_j-i-j+1]}.
\]
The  classical dimension is given by the same formula but with
integers in the place of quantum integers.